\documentclass{amsart}
\usepackage{graphicx}
\usepackage{amsmath}
\usepackage{amssymb}
\usepackage{amsthm}
\usepackage{esint}
\usepackage{mathtools}

\def\P{\mathbb{P}}
\def\S{\mathbb{S}}

\newtheorem{theorem}{Theorem}[section]
\newtheorem{lem}[theorem]{Lemma}

\theoremstyle{Corollary}

\newtheorem{prop}[theorem]{Proposition}

\theoremstyle{definition}
\newtheorem{definition}[theorem]{Definition}
\newtheorem{remark}[theorem]{Remark}

\numberwithin{equation}{section}

\begin{document}

\title{Convergence rate of the weighted Yamabe flow}

\author{Pak Tung Ho}
\address{Department of Mathematics, Sogang University, Seoul
04107, Korea}

\address{Korea Institute for Advanced Study, Seoul, 02455, Korea}

\email{paktungho@yahoo.com.hk}

\author{Jinwoo Shin}
\address{Korea Institute for Advanced Study, Hoegiro 85, Seoul 02455, Korea}
\email{shinjin@kias.re.kr}

\author{Zetian Yan}
\address{109 McAllister Building, Penn State University, University Park, PA 16802, USA}
\email{zxy5156@psu.edu}


\date{30th June, 2022.}


\begin{abstract}
The weighted Yamabe flow was the geometric flow introduced to study the weighted Yamabe
problem on smooth metric measure spaces.
Carlotto, Chodosh and Rubinstein have studied the convergence rate of the Yamabe flow.
Inspired by their result, we study in this paper the convergence rate of the weighted Yamabe flow.

\end{abstract}

\maketitle

\section{Introduction}

Given a closed (i.e. compact without boundary) Riemannian manifold
 $(M,g_0)$ of dimension $n\geq 3$, the \textit{Yamabe problem} is to find a metric conformal to $g_0$ such that
the scalar curvature $R_g$ of $g$ is constant.
This was solved by Aubin \cite{Aubin0}, Trudinger \cite{Trudinger} and Schoen \cite{Schoen}.

The \textit{Yamabe flow} is a geometric flow introduced to
study the Yamabe problem, which is defined as
\begin{equation}\label{YF}
\frac{\partial g(t)}{\partial t}=-(R_{g(t)}-r_{g(t)})g(t),
\end{equation}
where
\begin{equation*}
r_{g(t)}=\frac{\int_MR_{g(t)}dV_{g(t)}}{\int_MdV_{g(t)}}
\end{equation*}
is the average of the scalar curvature of $g(t)$.
The existence and convergence of the Yamabe flow
has been studied in \cite{Brendle4,Brendle5,Chow,Schwetlick&Struwe,Ye}.
See also  \cite{Azami&Razavi,Cheng&Zhu,Daneshvar&Razavi,Ho1,Ho2,Ho3,Ma&Cheng,Schulz} and references therein
for results related to the Yamabe flow.

In \cite{Chodosh}, Carlotto, Chodosh and Rubinstein
studied the rate of convergence of the Yamabe flow (\ref{YF}).
In particular, they proved the following:

\begin{theorem}[Theorem 1 in \cite{Chodosh}]\label{thm1}
Assume $g(t)$ is a solution of the Yamabe flow \eqref{YF}
that converges in $C^{2,\alpha}(M,g_\infty)$ to $g_\infty$ as $t\to\infty$ for some $\alpha\in (0,1)$.
Then there is a $\delta>0$ depending only on $g_\infty$ such that:\\
(i) If $g_\infty$ is an integrable critical point, then the convergence occurs at an exponential rate, that is
$$\|g(t)-g_\infty\|_{C^{2,\alpha}(M,g_\infty)}\leq Ce^{-\delta t}$$
for some constant $C>0$ depending on $g(0)$. \\
(ii) In general, the rate of convergence cannot be worse than polynomial, that is
$$\|g(t)-g_\infty\|_{C^{2,\alpha}(M,g_\infty)}\leq C(1+t)^{-\delta}$$
for some constant $C>0$ depending on $g(0)$.
\end{theorem}

\begin{theorem}[Theorem 2 in \cite{Chodosh}]\label{thm2}
Assume that $g_\infty$ is a nonintegrable critical point of the Yamabe energy with
order of integrability $p\geq 3$. If $g_\infty$ satisfies the Adams-Simon positive condition
$AS_p$, then there exists metric $g(0)$ conformal to $g_\infty$ such that
the solution $g(t)$ of the Yamabe flow \eqref{YF}
starting from $g(0)$ exists for all time and converges
in $C^\infty(M,g_\infty)$ to $g_\infty$
as $t\to\infty$. The convergence occurs ``slowly" in the sense that
$$C^{-1}(1+t)^{-\frac{1}{p-2}}\leq \|g(t)-g_\infty\|_{C^{2,\alpha}(M,g_\infty)}\leq C(1+t)^{-\frac{1}{p-2}}$$
for some constant $C>0$.
\end{theorem}

We refer the readers to \cite[Definition 8]{Chodosh}
and \cite[Definition 10]{Chodosh}
respectively for the precise definitions of integrable critical points
and Adams-Simon positive condition $AS_p$.

To explain our results requires some terminology.
A \textit{smooth metric measure space} is a four-tuple $(M^n, g, e^{-\phi} dV_g, m)$ of a Riemannian manifold $(M^n,g)$, a smooth measure $e^{-\phi}dV_g$ determined by a function $\phi\in C^{\infty}(M)$ and the Riemannian volume element of $g$, and a dimensional parameter $m\in [0,\infty]$. In the case $m=0$, we require $\phi=0$.

The \textit{weighted scalar curvature} of
a smooth metric measure space $(M,g,e^{-\phi}dV_g, m)$
is defined as
\begin{equation}\label{1.0}
R_\phi^m:=R_g+2\Delta_g\phi-\frac{m+1}{m}|\nabla_g\phi|^2_g,
\end{equation}
where $R_g$ is the scalar curvature of $g$,
$\Delta_g$ and $\nabla_g$ are respectively the Laplacian and the gradient of $g$.

Conformal equivalence between smooth metric measure spaces are defined as the following, see \cite{Case1} for more details.
\begin{definition}\label{condef}
Smooth metric measure spaces $(M^n, g, e^{-\phi}dV_g, m)$ \\and $(M^n, \hat{g}, e^{-\hat{\phi}}dV_{\hat{g}}, m)$ are conformally equivalent if there is a smooth function $\sigma\in C^{\infty}(M)$ such that
\begin{equation}\label{1.2}
(M^n, \hat{g}, e^{-\hat{\phi}}dV_{\hat{g}}, m)=(M^n, e^{\frac{2}{m+n-2}\sigma}g, e^{\frac{m+n}{m+n-2}\sigma}e^{-\phi}dV_g, m).
\end{equation}
In the case $m=0$, conformal equivalence is defined in the classical sense.
\end{definition}
If we denote $e^{\frac{1}{2}\sigma}$ by $w$, (\ref{1.2}) is equivalent to
\begin{equation}
(M^n, \hat{g}, e^{-\hat{\phi}}dV_{\hat{g}}, m)=(M^n,w^{\frac{4}{m+n-2}}g,w^{\frac{2(m+n)}{m+n-2}}e^{-\phi}dV_{g}, m),
\end{equation}
which is an alternative way to formulate the conformal equivalence of smooth metric measure spaces.

If $(M,g,e^{-\phi}dV_g, m)$
and $(M,g_0,e^{-\phi_0}dV_{g_0}, m)$ are conformal in the sense of (\ref{1.2}),
then their weighted scalar curvatures are related by (see (2.2) in \cite{Yan} for example)
\begin{equation}\label{1.6}
-\frac{4(n+m-1)}{n+m-2}\Delta_{\phi_0} w+R_{\phi_0}^mw=R_\phi^m w^{\frac{m+n+2}{m+n-2}},
\end{equation}
where
$$\Delta_{\phi_0}:=\Delta-\nabla \phi_0$$
is the \textit{weighted Laplacian} of $(M,g_0,e^{-\phi_0}dV_{g_0}, m)$, i.e.
$$\Delta_{\phi_0}\psi=\Delta_{g_0}\psi-\langle\nabla_{g_0} \phi_0,\nabla_{g_0}\psi\rangle~~\mbox{ for any }\psi\in C^\infty(M).$$
The fact that we will constantly use throughout this paper is that
the weighted Laplacian $\Delta_{\phi_0}$
is formally self-adjoint with respect to the measure $e^{-\phi_0}dV_{g_0}$
(c.f. \cite{Case1}).

Given a compact smooth metric measure space $(M,g_0,e^{-\phi_0}dV_{g_0}, m)$,
the \textit{weighted Yamabe problem} is to find
another smooth metric measure space $(M,g,e^{-\phi}dV_g, m)$
conformal to $(M,g_0,e^{-\phi_0}dV_{g_0}, m)$ such that
its weighted scalar curvature $R^m_\phi$ is constant.
The weighted Yamabe problem in this article is different from that introduced by Case in \cite{Case1}.
See also \cite{Case2,Case3,deSouza,Munoz} for more results related to Case's weighted Yamabe problem.

Similar to the Yamabe flow,
the \textit{weighted Yamabe flow}
is the geometric flow used to study the weighted Yamabe problem.
This was first introduced by Yan in \cite{Yan}.
More precisely, the weighted Yamabe flow is the evolution equation defined on
$(M,g(t),e^{-\phi(t)}dV_{g(t)}, m)$
given by
\begin{equation}\label{1.3}
\begin{split}
  \begin{dcases}
    \frac{\partial g}{\partial t} &=(r^m_{\phi}-R^m_{\phi})g, \\
\frac{\partial \phi}{\partial t} &=\frac{m}{2}(R^m_{\phi}-r^m_{\phi}),
  \end{dcases}
\end{split}
\end{equation}
where $r^m_{\phi}$ is the mean value of $R^m_{\phi}$; i.e.
\begin{equation}\label{1.4}
r^m_{\phi}=\frac{\int_M R^m_{\phi} e^{-\phi}dV_g}{\int_M e^{-\phi}dV_g}.
\end{equation}
In \cite{Yan}, Yan proved that the weighted Yamabe flow
(\ref{1.3}) exists for all time and converges to a metric
with constant weighted scalar curvature.

Inspired
by the results of Carlotto, Chodosh and Rubinstein about the convergence rate of Yamabe flow,
i.e. Theorems \ref{thm1} and \ref{thm2} mentioned above,
 we study in this paper the rate of convergence of the weighted Yamabe flow (\ref{1.3}).

The following theorems are the main results in this paper:

\begin{theorem}\label{main1}
Assume that $(g(t),\phi(t))$ is a solution to the weighted Yamabe flow
that is converging in $C^{2,\alpha}(M,g_\infty)$ to $(g_\infty,\phi_\infty)$
as $t\to\infty$ for some $\alpha\in (0,1)$. Then, there is $\delta>0$ depending only on $(g_\infty,\phi_\infty)$ such that\\
(i) If $(g_\infty,\phi_\infty)$ is an integrable critical point, then the convergence
occurs at an exponential rate
$$\|(g(t),\phi(t))-(g_\infty,\phi_\infty)\|_{C^{2,\alpha}(M,g_\infty)}\leq Ce^{-\delta t},$$
for some constant $C>0$ depending on $(g(0),\phi(0))$.\\
(ii) In general, the convergence cannot be worse than a polynomial rate
 $$\|(g(t),\phi(t))-(g_\infty,\phi_\infty)\|_{C^{2,\alpha}(M,g_\infty)}\leq C(1+t)^{-\delta},$$
for some constant $C>0$ depending on $(g(0),\phi(0))$, where
\begin{displaymath}
    \|(g(t),\phi(t))-(g_\infty,\phi_\infty)\|_{C^{2,\alpha}(M,g_\infty)}=\|g(t)-g_\infty\|_{C^{2,\alpha}(M,g_\infty)}+\|\phi(t)-\phi_\infty\|_{C^{2,\alpha}(M,g_\infty)}.
\end{displaymath}
\end{theorem}

\begin{theorem}\label{main2}
Assume that $(g_\infty,\phi_\infty)$ is a non-integrable critical point of
the energy functional $E$ with order of integrability $p\geq 3$. If $(g_\infty,\phi_\infty)$
satisfies the Adam-Simon positivity condition $AS_p$, then there exists
a metric-measure structure $(g(0),\phi(0))$ conformal to $(g_\infty,\phi_\infty)$
such that the weighted Yamabe flow
$(g(t),\phi(t))$ starting from $(g(0),\phi(0))$ exists for all time and converges in $C^\infty(M,g_\infty)$ to $(g_\infty,\phi_\infty)$
as $t\to\infty$. The convergence occurs ``slowly" in the sense that
$$C(1+t)^{-\frac{1}{p-2}}\leq \|(g(t),\phi(t))-(g_\infty,\phi_\infty)\|_{C^{2}(M,g_\infty)}\leq C(1+t)^{-\frac{1}{p-2}}$$
for some constant $C>0$.
\end{theorem}

The structure of this article is the following. Section \ref{section2} is devoted to fixing the notation and recalling some backgrounds about the normalized Yamabe functional, its analyticity and the Lyapunov--Schmidt reduction near a critical point. In particular, the precise definitions of integrable critical point and
the Adam-Simon positivity condition $AS_p$ can be found there. In Section \ref{section3}, we use the {\L}ojasiewicz-Simon inequality to prove Theorem \ref{thm1}. Next, in Section \ref{section4} we study polynomial convergence phenomena for nonintegrable critical points, and in Section \ref{section5}, we construct example of Riemannian manifolds
which satisfies the condition $AS_3$.
This allows us to conclude that
there exists a weighted Yamabe flow
converging  exactly at a polynomial rate described in Theorem \ref{main2}.

\section{Definitions and Preliminaries}\label{section2}\sloppy
\begin{definition}
On a smooth metric measure space $(M^n, g, e^{-\phi}dV_g, m)$, which is conformal to $(M^n, g_0, e^{-\phi_0}dV_{g_0}, m)$ in the sense of Definition \ref{condef},
\begin{displaymath}
(M^n, g, e^{-\phi}dV_g, m)=(M^n, w^{\frac{4}{n+m-2}}g_0, w^{\frac{2(m+n)}{n+m-2}}e^{-\phi_0}dV_{g_0}, m),
\end{displaymath}
analogous to the classical Yamabe problem, we define the normalized energy functional $E(w)$ as
\begin{equation}\label{energy}
     E_{(g_0,\phi_0)}(w)=\frac{\int_M \left(\frac{4(n+m-1)}{n+m-2}L^m_{\phi_0}w,w\right)e^{-\phi_0}dV_{g_0} }{\left( \int_M w^{\frac{2(n+m)}{n+m-2}}e^{-\phi_0}dV_{g_0}\right)^{\frac{n+m-2}{n+m}}},
\end{equation}
where $L^m_{\phi_0}$ is the weighted conformal Laplacian on $(M^n, g_0, e^{-\phi_0}dV_{g_0}, m)$
\begin{displaymath}
L^m_{\phi_0}=-\Delta_{\phi_0}+\frac{n+m-2}{4(n+m-1)}R^m_{\phi_0}.
\end{displaymath}
\end{definition}

\begin{remark}{The normalized total weighted scalar curvature:}
Under the setting of Definition \ref{energy}, by the transformation law of the weighted scalar curvature in (\ref{1.6}), the normalized energy $E_{(g_0,\phi_0)}(w)$ is exactly the normalized total weighted scalar curvature of $(M^n, g, e^{-\phi}dV_g, m)$; i.e.
\begin{displaymath}
E_{(g_0,\phi_0)}(w)=E_{(g,\phi)}(1)= \frac{\int_M R^m_{\phi}e^{-\phi}dV_{g} }{\left( {\rm Vol}\left(M^n, e^{-\phi}dV_g\right)\right)^{\frac{n+m-2}{n+m}}}.
\end{displaymath}
\end{remark}

Along the flow (\ref{1.3}), the volume
$\displaystyle\int_M e^{-\phi(t)}dV_{g(t)}$ is preserved.
Indeed,
it follows from (\ref{1.3}) and (\ref{1.4}) that
\begin{equation}\label{2.1}
\frac{d}{dt}\left(\int_Me^{-\phi(t)} dV_{g(t)}\right)
=\frac{n+m}{2}\int_M (r_{\phi(t)}^m-R_{\phi(t)}^m)dV_{g(t)}=0.
\end{equation}
Since the flow (\ref{1.3}) preserves the conformal structure,
we can write the solution as
\begin{equation}\label{2.4}
(M,g(t),e^{-\phi(t)}dV_{g(t)}, m)=(M,u(t)^{\frac{4}{m+n-2}}g_\infty,u(t)^{\frac{2(m+n)}{m+n-2}}e^{-\phi_\infty}dV_{g_\infty}, m).
\end{equation}
We remark that this implies that
\begin{equation}\label{2.4b}
  \phi(t)=\phi_\infty-\frac{2m}{m+n-2}\ln u(t).
\end{equation}

Therefore, we assume that the volume of $(M, g_\infty, e^{-\phi_\infty}dV_{g_\infty}, m)$
satisfying
\begin{equation}\label{2.5}
\int_Me^{-\phi_\infty}dV_{g_\infty}=1,
\end{equation}
then it follows from (\ref{2.1}) that
\begin{equation}\label{2.2}
\int_Me^{-\phi(t)}dV_{g(t)}=1~~\mbox{ for all }t\geq 0.
\end{equation}
In view of (\ref{2.4}) and \eqref{2.4b}, the weighted Yamabe flow (\ref{1.3}) reduces to the following evolution equation for the conformal factor:
\begin{equation}\label{1.8}
\frac{\partial}{\partial t}u(t)=\frac{n+m-2}{4} (r_{\phi(t)}^m-R_{\phi(t)}^m)u(t).
\end{equation}
Together this with (\ref{2.2}), we find
\begin{equation}\label{1.5}
\frac{d}{dt}r_{\phi(t)}^m=\frac{d}{dt}E(u(t))
=-\frac{n+m-2}{2} \int_M(R_{\phi(t)}^m-r_{\phi(t)}^m)^2e^{-\phi(t)}dV_{g(t)}\leqslant 0
\end{equation}
along the flow.

Consider the following unit volume conformal class associated to $(g_\infty,\phi_\infty)$:
\begin{equation*}
  \begin{split}
    [(g_\infty,\phi_\infty)]_1=\left\{(w^{\frac{4}{n+m-2}}g_\infty,\phi_\infty-\frac{2m}{m+n-2}\ln w): 0<w\in C^{2,\alpha}(M)\right.,\\
 \left.\int_M w^{\frac{2(n+m)}{n+m-2}}e^{-\phi_\infty}dV_{g_\infty}=1\right\}.
  \end{split}
\end{equation*}
$$$$

In order to avoid ambiguities, we define the following notion: for $k\in\mathbb{N}$, we denote the $k$-th differential of the energy functional $E$ on $[(g_\infty,\phi_\infty)]_1$
at the point $w$ in the directions $v_1,..., v_k$ by
$$D^kE(w)[v_1,..., v_k].$$
As we will see below, the functional $v\mapsto D^kE(w)[v_1,..., v_{k-1},v]$
is in the image of $L^2(M,e^{-\phi_\infty}dV_{g_\infty})$ under the natural embedding onto $C^{2,\alpha}(M,g_\infty)'$.
Therefore, we will also write
$$D^kE(w)[v_1,..., v_{k-1}]$$
for this element of $L^2(M,e^{-\phi_\infty}dV_{g_\infty})$.
When $k=1$, we will drop the (second) brackets, and thus consider
$DE(w)\in L^2(M,e^{-\phi_\infty}dV_{g_\infty})$.

We may write the differential of $E$ restricted to $[(g_\infty,\phi_\infty)]_1$ as
\begin{equation}\label{diffE}
\begin{split}
\frac{1}{2}DE(w)[v]
&=\frac{1}{2}\left.\frac{d}{dt}E(w+tv)\right|_{t=0}\\
&=\frac{\int_M\big(\frac{4(n+m-1)}{n+m-2}\langle\nabla_{g_\infty}w,\nabla_{g_\infty}v\rangle+R_{\phi_\infty}^mwv\big)e^{-\phi_\infty}dV_{g_\infty}}
{\big(\int_Mw^{\frac{2(n+m)}{n+m-2}}e^{-\phi_\infty}dV_{g_\infty}\big)^{\frac{n+m-2}{n+m}}}\\
&\hspace{4mm}
-\frac{E(w)}{\int_Mw^{\frac{2(n+m)}{n+m-2}}e^{-\phi_\infty}dV_{g_\infty}}
\int_Mw^{\frac{n+m+2}{n+m-2}}ve^{-\phi_\infty}dV_{g_\infty}\\
&=\int_M\left(-\frac{4(n+m-1)}{n+m-2}\Delta_{\phi_\infty}w+R^m_{\phi_\infty}w-r^m_{\phi}w^{\frac{n+m+2}{n+m-2}}\right)ve^{-\phi_\infty}dV_{g_\infty}\\
&=\int_M\left(R^m_{\phi}-r^m_{\phi}\right)w^{\frac{n+m+2}{n+m-2}}ve^{-\phi_\infty}dV_{g_\infty}
\end{split}
\end{equation}
where
$$(M,g,e^{-\phi}dV_g, m)=(M,w^{\frac{4}{m+n-2}}g_\infty,w^{\frac{2(m+n)}{m+n-2}}e^{-\phi_\infty}dV_{g_\infty}, m).
$$
Thus, a unit volume metric-measure structure $(g,\phi)$ is a critical point for the energy
$E$ restricted to $[(g_\infty,\phi_\infty]_1$ exactly when
$(g,\phi)$ has constant weighted scalar curvature.

We now fix $(g_\infty,\phi_\infty)$ such that
(\ref{2.5}) holds and $(g_\infty,\phi_\infty)$ has constant weighted scalar curvature.
We denote by $\mathcal{CWSC}_1$ the set of unit volume constant weighted scalar curvature metric-measure structures
in $[(g_\infty,\phi_\infty)]_1$.
If we define the \textit{linearized weighted Yamabe operator} at $(g_\infty,\phi_\infty)$, $\mathcal{L}_\infty$,
by means of the formula
\begin{equation*}
  \begin{split}
    -\frac{4}{n+m-2}\int_M w\mathcal{L}_\infty v e^{-\phi_\infty}dV_{g_\infty}
:=&\frac{1}{2}D^2E(g_\infty,\phi_\infty)[v,w]\\
=&\frac{1}{2}\left.\frac{d}{dt}\left(DE(1+tw)[v]\right)\right|_{t=0}
  \end{split}
\end{equation*}
for $v\in C^2(M)$. A direct computation (see the Appendix) shows that
$$\mathcal{L}_\infty v=(n+m-1)\Delta_{\phi_\infty}v+R_{\phi_\infty}^mv.$$
We define $\Lambda_0:=\ker\mathcal{L}_\infty\subset L^2(M,e^{-\phi_\infty}dV_{g_\infty})$.

It follows from a classical theorem of spectral theory that $\Lambda_0$ is finite dimensional,
since it is the  eigenspace of the weighted Laplacian $\Delta_{\phi_\infty}$
for the eigenvalue  $\displaystyle\frac{R_{\phi_\infty}^m}{n+m-1}$.
We will write $\Lambda_0^{\perp}$
for the $L^2(M,e^{-\phi_\infty}dV_{g_\infty})$-orthogonal complement.

\medskip

It is crucial throughout this work that the functional $E$ is an analytic map in the sense of  \cite[Definition 8.8]{Zeidler}. More precisely, one can easily prove the following by expanding the denominator of $E$ in a power series: fix a metric-measure structure $(g_\infty,\phi_\infty)$ then  the functional $E$ is an analytic functional on $\{u\in C^{2,\alpha}(M,g_\infty):u>0\}$ in the sense that for each $w_0\in C^{2,\alpha}(M,g_\infty)$ with $w_0>0$, there is an $\epsilon>0$ and bounded multilinear operators
  \begin{equation*}
    E^{(k)}:C^{2,\alpha}(M,g_\infty)^{\times k}\rightarrow \mathbb{R}\textrm{ for each }k\geq 0
  \end{equation*}
  such that if $\|w-w_0\|_{C^{2,\alpha}}<\epsilon$, then $\sum_{k=0}^\infty\|E^{(k)}\|\cdot\|w-w_0\|^k_{C^{2,\alpha}}<\infty$ and
  \begin{equation*}
    E(w)=\sum_{k=0}^\infty E^{(k)}(\underbrace{w-w_0,\cdots,w-w_0}_{k\textrm{-times}})\textrm{ in }C^{2,\alpha}(M,g_\infty).
  \end{equation*}

\medskip
We need the following proposition from \cite[Section 3]{Simon}, which can be established with the help of the implicit function theorem:

\begin{prop}\label{prop7}
There is $\epsilon>0$ and an analytic map $$\Phi:\Lambda_0\cap \{v: \|v\|_{L^2}<\epsilon\}
\to C^{2,\alpha}(M,g_\infty)\cap \Lambda_0^\perp$$
such that $\Phi(0)=0$, $D\Phi(0)=0$,
\begin{equation}\label{eq2.7}
\sup_{\substack{
\|v\|_{L^2}<\epsilon,\\
\|w\|_{L^2}\leq 1}}\|D\Phi(v)[w]\|_{L^2}<1,
\end{equation}
and so that defining $\Psi(v)=1+v+\Phi(v)$,
we have that $\Psi(v)>0$, $\displaystyle\int_M \Psi(v)^{\frac{2(n+m)}{n+m-2}}e^{-\phi_\infty}dV_{g_\infty}=1$ and
$$\mbox{\emph{proj}}_{\Lambda_0^\perp}[DE(\Psi(v))]=
\mbox{\emph{proj}}_{\Lambda_0^\perp}
\left[\Big(R_{\phi}^m-r_{\phi}^m\Big)\Psi(v)^{\frac{n+m+2}{n+m-2}}\right]=0$$
where
\begin{equation*}
  (g,\phi)=(\Psi(v)^{\frac{4}{n+m-2}}g_\infty,\phi_\infty-\frac{2m}{m+n-2}\ln \Psi(v)).
\end{equation*}

Furthermore
$$\mbox{\emph{proj}}_{\Lambda_0}[DE(\Psi(v))]=
\mbox{\emph{proj}}_{\Lambda_0}
\left[\Big(R_{\phi}^m-r_{\phi}^m\Big)\Psi(v)^{\frac{n+m+2}{n+m-2}}\right]=DF(v),$$
where $F:\Lambda_0\cap\{v: \|v\|_{L^2}\leq \epsilon\}\to\mathbb{R}$ is defined by
$F(v)=E(\Psi(v))$. Finally, the intersection of $\mathcal{CWSC}_1$ with a small $C^{2,\alpha}(M,g_\infty)$-neighborhood
of $1$ coincides with
$$\mathcal{S}_0:=\{\Psi(v): v\in\Lambda_0, \|v\|_{L^2}<\epsilon, DF(v)=0\},$$
which is a real analytic subvariety (possible singular) of the following
$(\dim\Lambda_0)$-dimensional real analytic submanifold of $C^{2,\alpha}(M,g_\infty)$:
$$\mathcal{S}:\{\Psi(v): v\in\Lambda_0, \|v\|_{L^2}<\epsilon\}.$$
\end{prop}

We will refer to $\mathcal{S}$ as the \textit{natural constraint} for the problem.

\begin{definition}\label{def8}
  For $(g_\infty,\phi_\infty)\in \mathcal{CWSC}_1$, we say that $(g_\infty,\phi_\infty)$ is \textit{integrable} if for all $v\in\Lambda_0$, there is a path $w(t)\in C^2((-\epsilon,\epsilon)\times M,g_\infty)$ such that $(w(t)^\frac{4}{n+m-2}g_\infty,\phi_\infty-\frac{2m}{m+n-2}\ln w(t))\in \mathcal{CWSC}_1$ and $w(0)=1$, $w'(0)=v$. Equivalently, $(g_\infty,\phi_\infty)$ is integrable if and only if $\mathcal{CWSC}_1$ agrees with $\mathcal{S}$ in a small neighborhood of $1$ in $C^{2,\alpha}(M,g_\infty)$.
\end{definition}

We remark that the integrability defined in Definition \ref{def8} is equivalent to the functional $F$ (as defined in Proposition \ref{prop7}) being constant in a neighborhood of $0$ inside $\Lambda_0$ \cite[Lemma 1]{Adams&Simon}.

\begin{definition}
  If $\Lambda_0=\{0\}$, i.e. if $\mathcal{L}_\infty$ is injective, then we call $(g_\infty,\phi_\infty)$ a \textit{nondegenerate critical point}. On the other hand, if $\Lambda_0$ is nonempty, we call $(g_\infty,\phi_\infty)$ \textit{degenerate}.
\end{definition}

Note that if $(g_\infty,\phi_\infty)$ is a nondegenerate critical point, then $(g_\infty,\phi_\infty)$ is automatically integrable in the above sense.

\medskip

Now suppose that $(g_\infty,\phi_\infty)$ is a nonintegrable critical point. Because $F(v)=E(\Psi(v))$, defined in Proposition \ref{prop7},  is analytic, we may expand it in a power series
\begin{equation*}
  F(v)=F(0)+\sum_{j\geq p}F_j(v)
\end{equation*}
where $F_j$ is a degree-$j$ homogeneous polynomial on $\Lambda_0$ and $p$ is chosen so that $F_p$ is nonzero. We will call $p$ the \textit{order of integrability} of $(g_\infty,\phi_\infty)$. We will also need a further hypothesis for nonintegrable critical points introduced in \cite{Adams&Simon}.

\begin{definition}\label{def10}
  We say that $(g_\infty,\phi_\infty)$ satisfies the \textit{Adams-Simon positivity condition}, $AS_p$ for short (here $p$ is the order of integrability of $g_\infty$), if it is nonintegrable and $F_p|_{\mathbb{S}^k}$ attains a positive maximum for some $\hat{v}\in \mathbb{S}^k\subset \Lambda_0$. Recall that $F_p$ is the lowest-degree nonconstant term in the power series expansion of $F(v)$ around $0$ and $\mathbb{S}^k$ is the unit sphere with respect to the $L^2(M,e^{-\phi_\infty}dV_{g_\infty})$-norm in $\Lambda_0$.
\end{definition}

\medskip

An important observation is that when the order of integrability $p$ is odd, the Adams-Simon positivity condition  is always satisfied. Moreover the order of integrability (at a critical point of $E$) always satisfies $p\geq 3$. Furthermore, we will show in the Appendix that
\begin{equation}\label{1.10}
F_3(v)=-\frac{8(n+m+2)}{(n+m-2)^2}R_{\phi_\infty}^m\int_M v^3 e^{-\phi_\infty}dV_{g_\infty}.
\end{equation}

\medskip

\section{The rate of convergence}\label{section3}

One of the tools for controlling the rate of convergence of weighted Yamabe flow will be the {\L}ojasiewicz-Simon inequality stated in \cite[Definition 11]{Chodosh}.

\begin{prop}\label{prop13}
Suppose that $(g_\infty,\phi_\infty)$ satisfies \eqref{2.5}
and has constant weighted scalar curvature. There are $\theta\in(0,\frac{1}{2}]$,
$\epsilon>0$ and $C>0$ (both depending only on $n$ and $(g_\infty,\phi_\infty)$)
such that for $u\in C^{2,\alpha}(M,g_\infty)$ with
$\|u-1\|_{C^{2,\alpha}(M,g_\infty)}<\epsilon$, then
$$\Big|r_{\phi}^m-r_{\phi_\infty}^m\Big|^{1-\theta}\leq C\|DE(g,\phi)\|_{L^2(M,e^{-\phi_\infty}dV_{g_\infty})}$$
where
\begin{equation*}
  (g,\phi)=(u^\frac{4}{n+m-2}g_\infty,\phi_\infty-\frac{2m}{m+n-2}\ln u).
\end{equation*}

If $(g_\infty,\phi_\infty)$ is an integrable critical point, then $\theta=\frac{1}{2}$.
If $(g_\infty,\phi_\infty)$ is non-integrable,
then this holds for some $\theta\in(0,\frac{1}{p}]$, where $p$ is the order of integrability of
$(g_\infty,\phi_\infty)$.
\end{prop}
\begin{proof}
  To verify this, we will show that hypotheses of Proposition 12 in \cite{Chodosh} are satisfied for the energy functional $E$. We work with the Banach spaces ${\mathcal{B}}:=C^{2,\alpha}(M,g_\infty)$ and ${\mathcal{W}}:=L^2(M,e^{-\phi_\infty}dV_{g_\infty})$, and fix $U$ a small enough ball around $1$ in $C^{2,\alpha}(M,g_\infty)$ so that Proposition \ref{prop7} is applicable in $U$.

  Hypothesis (A) says that $\Lambda_0=\ker\mathcal{L}_\infty$ is complemented in $C^{2,\alpha}(M,g_\infty)$, which is immediate by the following argument. It's not hard to check that the $L^2$ projection map ${\rm{proj}}_{\Lambda_0}$ restricts to a continuous map from $C^{2,\alpha}(M,g_\infty)$ onto $\Lambda_0$ (since, of course, $C^{2,\alpha}(M,g_\infty) \hookrightarrow L^2(M,e^{-\phi_\infty}dV_{g_\infty})$ is a continuous embedding); from this, it follows that $\Lambda_0^{'}$ is complemented (by the map ${\rm{proj}}^{'}_{\Lambda_0}$) in the dual space $C^{2,\alpha}(M,g_\infty)^{'}$ may be canonically identified with $(\Lambda_0^{\perp})^{'}$.

  Hypothesis (B) is satisfied as follows: Consider the map
  \begin{equation}
      {\mathcal{W}}:=L^2(M,e^{-\phi_\infty}dV_{g_\infty}) \hookrightarrow  C^{2,\alpha}(M,g_\infty)^{'}, \quad f \mapsto \left(\psi\mapsto \int_M f\psi e^{-\phi_\infty}dV_{g_\infty}\right).
  \end{equation}
  (B1) This map is continuous.\\
  (B2) The map ${\rm{proj}}^{'}_{\Lambda_0}\in {\mathcal{B}}(C^{2,\alpha}(M,g_\infty)^{'})$ leaves $L^2(M,e^{-\phi_\infty}dV_{g_\infty})$ invariant; here we are considering the composition
  \begin{equation}
      {\rm{proj}}_{\Lambda_0}: C^{2,\alpha}(M,g_\infty)\to \Lambda_0 \hookrightarrow C^{2,\alpha}(M,g_\infty).
  \end{equation}
  (B3) The fact that $DE\in C^1(U,L^2(M,e^{-\phi_\infty}dV_{g_\infty}))$ follows from the explicit form of $DE$ given above.\\
  (B4) Finally, we have to verify that range $\mathcal{L}_\infty=(\Lambda_0^{\perp})^{'}\cap L^2(M,e^{-\phi_\infty}dV_{g_\infty})$. The fact that range $\mathcal{L}_\infty \subset (\Lambda_0^{\perp})^{'}\cap L^2(M,e^{-\phi_\infty}dV_{g_\infty})$ is obvious because $\mathcal{L}_\infty$ is formally self-adjoint on $L^2(M,e^{-\phi_\infty}dV_{g_\infty})$.  The other inclusion follows from the $L^2$ spectral decomposition of $\mathcal{L}_\infty$.

  Therefore, to prove the {\L}ojasiewicz-Simon inequality with exponent $\theta\in (0,\frac{1}{2}]$, it suffices to check hypotheses (C), i.e. that the energy functional $E$ restricted to the natural constraint satisfies the {\L}ojasiewicz-Simon inequality  with exponent $\theta\in (0,\frac{1}{2}]$. Recall that in Proposition \ref{prop7} we have defined $F(v)=E(\Psi(v))$. In the integrable case, clearly, $F(v)= F(0) $, so $F$ satisfies the {\L}ojasiewicz-Simon inequality  with exponent $\frac{1}{2}$.

  In general, by definition, $F$ is an analytic function whose power series has its first nonzero term of degree $p$. Similar to \cite[Proposition 13]{Chodosh}, we may conclude that $F$ satisfies the {\L}ojasiewicz-Simon inequality  with exponent $\frac{1}{p}$.

  The claim follows from the fact that $E(u)=r^m_{\phi}$ under the volume normalization.
\end{proof}

Now we show how the {\L}ojasiewicz-Simon inequality yields quantitative estimates on the rate of convergence of the weghted Yamabe flow.

\begin{proof}[Proof of Theorem \ref{main1}]
We consider the weighted Yamabe flow
 $(g(t),\phi(t))=(u(t)^{\frac{4}{n+m-2}}g_\infty,\phi_\infty-\frac{2m}{m+n-2}\ln u(t))$
which converges to $(g_\infty,\phi_\infty)$ in $C^{2,\alpha}(M,g_\infty)$ as $t\to\infty$.
In Proposition \ref{prop13}, we have shown that there is a {\L}ojasiewicz-Simon
inequality near $(g_\infty,\phi_\infty)$ for some $\theta\in (0,\frac{1}{2}]$.
We emphasize that if we regard $DE(u(t))$ as an element of $L^2(M,e^{-\phi_\infty}dV_{g_\infty})$, then
\begin{equation}\label{1.9}
DE(u(t))=2\big(R_{\phi(t)}^m-r_{\phi(t)}^m\big)u(t)^{\frac{n+m+2}{n+m-2}}.
\end{equation}
Choose $t_0$ large enough to apply the {\L}ojasiewicz-Simon inequality. In other words,
 so that $\|u(t)-1\|_{C^{2,\alpha}(M,g_\infty)}\leq \epsilon$ for all $t\geq t_0$.
This together with (\ref{1.5}) and Proposition \ref{prop13} implies that
\begin{equation}\label{2.6}
\begin{split}
\frac{d}{dt}\big(r_{\phi(t)}^m-r_{\phi_\infty}^m\big)
&=-\frac{n+m-2}{2}\int_M(R_{\phi(t)}^m-r_{\phi(t)}^m)^2u(t)^{\frac{2(n+m)}{n+m-2}}e^{-\phi_\infty}dV_{g_\infty}\\
&\leq -c\int_M(R_{\phi(t)}^m-r_{\phi(t)}^m)^2u(t)^{\frac{2(n+m+2)}{n+m-2}}e^{-\phi_\infty}dV_{g_\infty}\\
&=-c\big\|DE(u(t))\big\|_{L^2(M,e^{-\phi_\infty}dV_{g_\infty})}^2\\
&\leq -c\big(r_{\phi(t)}^m-r_{\phi_\infty}^m\big)^{2-2\theta},
\end{split}
\end{equation}
where $c>0$ is a constant depending only on $n$ and $(g_\infty,\phi_\infty)$ (that we let change from line to line).
Let us first assume that the {\L}ojasiewicz-Simon inequality is satisfied with $\theta=\frac{1}{2}$, i.e.
we are in the integrable case. Then (\ref{2.6}) yields
$0\leq r_{\phi(t)}^m-r_{\phi_\infty}^m\leq C e^{-2\delta t}$, for some $\delta>0$ depending only on $n$ and $(g_\infty,\phi_\infty)$,
and $C>0$ depending on $(g(0),\phi(0))$ (chosen so that this actually holds for all $t\geq 0$).
On the other hand, if {\L}ojasiewicz-Simon inequality holds with $\theta\in (0,\frac{1}{2})$,
then the same argument shows that
$r_{\phi(t)}^m-r_{\phi_\infty}^m\leq C(1+t)^{\frac{1}{2\theta-1}}$.

Exploiting the fact that the flow converges in $C^2$, we may use the
{\L}ojasiewicz-Simon inequality to compute
\begin{equation*}
\begin{split}
\frac{d}{dt}\big(r_{\phi(t)}^m-r_{\phi_\infty}^m\big)^\theta
&=\theta\big(r_{\phi(t)}^m-r_{\phi_\infty}^m\big)^{\theta-1}
\frac{d}{dt}\big(r_{\phi(t)}^m-r_{\phi_\infty}^m\big)\\
&\leq -c\,\theta\big(r_{\phi(t)}^m-r_{\phi_\infty}^m\big)^{\theta-1}\big\|DE(u(t))\big\|_{L^2(M,e^{-\phi_\infty}dV_{g_\infty})}^2\\
&\leq -c\,\theta\|DE(u(t))\big\|_{L^2(M,e^{-\phi_\infty}dV_{g_\infty})}\\
&\leq -c\,\theta\left\|\frac{\partial u(t)}{\partial t}\right\|_{L^2(M,e^{-\phi_\infty}dV_{g_\infty})},
\end{split}
\end{equation*}
where we have used  (\ref{1.8}) and (\ref{1.9}) in the last equality.
Thus, if $\theta=\frac{1}{2}$ (recall $\lim_{t\to\infty}u(t)=1$), then
\begin{equation*}
\begin{split}
\|u(t)-1\|_{L^2(M,e^{-\phi_\infty}dV_{g_\infty})} & \leq\int_t^\infty\left\|\frac{\partial u(s)}{\partial s}\right\|_{L^2(M,e^{-\phi_\infty}dV_{g_\infty})}ds \\
     & \leq -c\int_t^\infty\frac{d}{ds}\left[\big(r_{\phi(s)}^m-r_{\phi_\infty}^m\big)^{\frac{1}{2}}\right]ds\\
     &=c\big(r_{\phi(t)}^m-r_{\phi_\infty}^m\big)^{\frac{1}{2}}\leq Ce^{-\delta t}.
\end{split}
\end{equation*}
If $\theta\in (0,\frac{1}{2})$, a similar computation yields $\|u(t)-1\|_{L^2(M,e^{-\phi_\infty}dV_{g_\infty})}\leq C(1+t)^{-\frac{\theta}{1-2\theta}}$.

To obtain $C^2$ estimates, we may interpolate between $L^2(M,e^{-\phi}g)$ and $W^{k,2}(M,e^{-\phi}g)$
for $k$ large enough:
interpolation \cite[Theorem 6.4.5]{Bergh} and Sobolev embedding yields
some constant $\eta\in (0,1)$ so that
$$\|u(t)-1\|_{C^{2,\alpha}(M,e^{-\phi_\infty}dV_{g_\infty})}\leq \|u(t)-1\|_{L^2(M,e^{-\phi_\infty}dV_{g_\infty})}^\eta\|u(t)-1\|_{W^{k,2}(M,e^{-\phi_\infty}dV_{g_\infty})}^{1-\eta}.$$
Because $u(t)$ converges to $1$ in $C^{2,\alpha}$
(and thus in $C^\infty$ by parabolic Schauder estimates and bootstrapping),
the second term is uniformly bounded.
Thus, exponential (polynomial) decay of the $L^2$ norm gives
exponential (polynomial) decay of the $C^{2,\alpha}$ norm as well.

Since $(g(t), \phi(t))=(u(t)^\frac{4}{n+m-2}g_\infty,\phi_\infty-\frac{2m}{m+n-2}\ln u(t))$, we immediately have
\begin{equation}
    \|(g(t),\phi(t))-(g_\infty,\phi_\infty)\|_{C^{2,\alpha}(M,g_\infty)}\leq Ce^{-\delta t},
\end{equation}
for some constant $C>0$ depending on $(g(0),\phi(0))$.
\end{proof}

\section{Slowly converging weighted Yamabe flow}\label{section4}

In this section, we show that, given a nonintegrable critical point $(g_\infty,\phi_\infty)$ satisfying a particular hypothesis, there exists a weighted Yamabe flow $(g(t),\phi(t))$ such that $(g(t),\phi(t))$ converges to $(g_\infty,\phi_\infty)$ exactly  at a polynomial rate.

This section is organized as follows: In section \ref{sec3.1}, we show that the weighted Yamabe flow can be represented by two different flows. To be more specific, we will project the flow equation to the kernel $\Lambda_0$ of $\mathcal{L}_\infty$ and its orthogonal complement $\Lambda_0^\perp$, respectively. In section \ref{sec3.2} and section \ref{sec3.3}, we solve the kernel-projected flow and the kernel-orthogonal projected flow, respectively. In section \ref{sec3.4}, we combine all the previous results to prove Theorem \ref{main2}.

\subsection{Projecting the weighted Yamabe flow with estimates}\label{sec3.1}

Here and in the sequel we will always use $f'(t)$ to denote the time derivative of a function $f(t)$. We will skip the proof of the following lemma, for its proof is the same as that of \cite[Lemma 15]{Chodosh}.

\begin{lem}\label{lem15}
  Assume that $(g_\infty,\phi_\infty)$ satisfies $AS_p$ as defined in Definition \ref{def10}, i.e. $F_p|_{\mathbb{S}^k}$ achieves a positive maximum for some point $\hat{v}$ in the unit sphere $\mathbb{S}^k\subset \Lambda_0$. Then, for any fixed $T\geq 0$, the function
  \begin{equation}\label{eq6}
    \varphi(t):=\varphi(t,T)=(T+t)^{-\frac{1}{p-2}}\left(\frac{8}{(n+m-2)p(p-2)F_p(\hat{v})}\right)^\frac{1}{p-2}\hat{v}
  \end{equation}
  solves $\frac{8}{n+m-2}\varphi'+DF_p(\varphi)=0$.
\end{lem}
We define the parabolic $C^{k,\alpha}$ norm on $(t,t+1)\times M$ as follows: for $\alpha\in(0,1)$, we define the seminorm
\begin{equation*}
  |f(t)|_{C^{0,\alpha}}:=\sup_{\substack{(s_i,x_i)\in(t,t+1)\times M \\ (s_1,x_1)\neq(s_2,x_2)}}\frac{|f(s_1,x_1)-f(s_2,x_2)|}{(d_{g_\infty}(x_1,x_2)^2+|t_1-t_2|)^\frac{\alpha}{2}}
\end{equation*}
and for $k\geq 0$ and $\alpha\in(0,1)$, we define the norm
\begin{equation}\label{eq8}
  \|f(t)\|_{C^{k,\alpha}}:=\sum_{|\beta|+2j\leq k}\sup_{(t,t+1)\times M}|D^\beta_x D^j_t f|+\sum_{|\beta|+2j=k}|D^\beta_x D^j_t f|_{C^{0,\alpha}}
\end{equation}
where the norm and derivatives in the sum are taken with respect to $g_\infty$.
When we mean an alternative norm, we will always indicate the domain.

\begin{lem}\label{appdA}
For the functional $E$, for $w$ such that $\|w-1\|_{C^{2,\alpha}}<1$, there holds
\begin{equation}\label{eq27}
  \|D^3E(w)[u,v]\|_{C^{0,\alpha}}\leq C\|u\|_{C^{2,\alpha}}\|v\|_{C^{2,\alpha}}
\end{equation}
for some uniform constant $C>0$. Furthermore, for $w_1,w_2$ such that $\|w_i-1\|_{C^{2,\alpha}}<1$, we have
\begin{equation*}
\begin{split}
  \|D^3E(w_1)[v,v]-D^3E(w_2)[u,u]\|_{C^{0,\alpha}}\leq& C(\|w_1\|_{C^{2,\alpha}}+\|w_2\|_{C^{2,\alpha}})\\
&\qquad \times (\|u\|_{C^{2,\alpha}}+\|v\|_{C^{2,\alpha}})\|u-v\|_{C^{2,\alpha}}
\end{split}
\end{equation*}
for some uniform constant $C>0$.
\end{lem}
\begin{proof}
It follows from the following computation for  $D^3E$ which proved in the Appendix:
\begin{equation*}
  D^3E(1)[v,u,z]=-\frac{8(n+m+2)}{(n+m-2)^2}R_{\phi_\infty}^m\int_M vuze^{-\phi_\infty}dV_{g_\infty}.
\end{equation*}
\end{proof}
\begin{lem}\label{lemma16}
  There exists $T_0>0$, $\epsilon_0>0$ and $c>0$, all depending on $(g_\infty,\phi_\infty)$ and $\hat{v}$, such that the following holds: Fix $T>T_0$. Then, for $\varphi(t)$ as in Lemma \ref{lem15} and $w\in C^{2,\alpha}(M\times[0,\infty))$, and $u:=1+\varphi+w^\top+\Phi(\varphi+w^\top)+w^\perp$ where $w^\top:=\textrm{proj}_{\Lambda_0}(w)$ and $w^\perp:=\textrm{proj}_{\Lambda_0^\perp}(w)$, the function
\begin{equation*}
  E_0^\top(w):=\mbox{\emph{proj}}_{\Lambda_0}\left[DE(u)u^{-\frac{4}{n+m-2}}-DE(u)\right]
\end{equation*}
satisfies
\begin{equation*}
  \begin{split}
  &\left\|E_0^\top(w)\right\|_{C^{0,\alpha}}\leq c\left\{(T+t)^{-\frac{p-1}{p-2}}+\|w^\top\|^{p-1}_{C^{0,\alpha}}+\|w^\perp\|_{C^{2,\alpha}}\right\}\left\{(T+t)^{-\frac{1}{p-2}}+\|w\|_{C^{2,\alpha}}\right\},\\
&\left\|E^\top_0(w_1)-E_0^\top(w_2)\right\|_{C^{0,\alpha}}\\
&\qquad\leq c\left\{(T+t)^{-\frac{p-1}{p-2}}+\|w_1^\top\|^{p-1}_{C^{0,\alpha}}+\|w_2^\top\|^{p-1}_{C^{0,\alpha}}+\|w_1^\perp\|_{C^{2,\alpha}}+\|w_2^\perp\|_{C^{2,\alpha}}\right\}\\
&\qquad \qquad \times \|w_1-w_2\|_{C^{2,\alpha}}\\
&+c\left\{(T+t)^{-\frac{1}{p-2}}+\|w_1\|_{C^{2,\alpha}}+\|w_2\|_{C^{2,\alpha}}\right\}\left(\|w_1^\top\|^{p-2}_{C^{0,\alpha}}+\|w_2^\top\|^{p-2}_{C^{0,\alpha}}\right)\\
&\qquad\qquad \times \|w_1^\top-w_2^\top\|_{C^{0,\alpha}}\\
&+c\left\{(T+t)^{-\frac{1}{p-2}}+\|w_1\|_{C^{2,\alpha}}+\|w_2\|_{C^{2,\alpha}}\right\}\|w_1^\perp-w_2^\perp\|_{C^{2,\alpha}}.
\end{split}
\end{equation*}
Identical estimates hold for $E_0^\perp(w):=\mbox{\emph{proj}}_{\Lambda_0^\perp}\left[DE(u)u^{-\frac{4}{n+m-2}}-DE(u)\right]$. Here, we are using the parabolic H\"{o}lder norms on $(t,t+1)\times M$ as defined above; the bounds hold for each fixed $t\geq 0$, with the constants independent of $T$ and $t$.
\end{lem}
\begin{proof}
Let $\eta:=\varphi+w^\top+\Phi(\varphi+w^\top)+w^\perp$. Then one can easily see that  $u=1+\eta$ and
\begin{equation*}
  \frac{d}{ds}(1+s\eta)^{-\frac{4}{n+m-2}}=-\frac{4}{n+m-2}(1+s\eta)^{-\frac{n+m+2}{n+m-2}}\eta.
\end{equation*}

Thus we have
\begin{equation*}
  \begin{split}
  u^{-\frac{4}{n+m-2}}&=1-\frac{4}{n+m-2}\int_0^1(1+s\eta)^{-\frac{n+m+2}{n+m-2}}\eta ds.
\end{split}
\end{equation*}
So, letting $E_0(w):=DE(u)u^{-\frac{4}{n+m-2}}-DE(u)$, we have
\begin{equation}\label{eq9}
  \begin{split}
  \|E_0(w)\|_{C^{0,\alpha}}=& c \left\|DE(u)\int_0^1(1+s\eta)^{-\frac{n+m+2}{n+m-2}}\eta ds\right\|_{C^{0,\alpha}}\\
  \leq & c\|DE(u)\|_{C^{0,\alpha}}\left( \|\varphi\|_{C^{0,\alpha}}+\|w^{\perp}\|_{C^{0,\alpha}}+\|\Phi(\varphi+w^\top)\|_{C^{0,\alpha}}+\|w^\top\|_{C^{0,\alpha}}\right)\\
  \leq & c\|DE(u)\|_{C^{0,\alpha}}\left((T+t)^{-\frac{1}{p-2}}+\|w^\top\|_{C^{0,\alpha}}+\|w^\perp\|_{C^{0,\alpha}}\right),
\end{split}
\end{equation}
where we have used the fact that $\Phi(0)=0$ and $\Phi$ is an analytic map.
It follows from  Taylor's theorem and Proposition \ref{prop7} that, for $\psi_{s,r}=1+r\left[\varphi+w^\top+\Phi(\varphi+w^\top)+sw^\perp\right]$,
\begin{equation}\label{taylor1}
  \begin{split}
  DE(u)=&DE(\Psi(\varphi+w^\top))+\int_0^1D^2E(\psi_{s,1})[w^\top]ds\\
  =&DF(\varphi+w^\top)-\frac{8}{n+m-2}\mathcal{L}_\infty w^\perp\\
                   &\qquad    +\int_0^1\int_0^sD^3E(\psi_{s,\tilde{s}})[w^\perp,\varphi+w^\top+\Phi(\varphi+w^\top)+sw^\perp]d\tilde{s}ds.
\end{split}
\end{equation}
 Now, observe that $DF(0)=D^2F(0)=\cdots D^{p-1}F(0)=0$, where $p$ is the order of integrability. Therefore, by Taylor's theorem, we have
\begin{equation}\label{taylor2}
  \|DF(\varphi+w^\top)\|_{C^{0,\alpha}}\leq c\|\varphi+w^\top\|^{p-1}_{C^{0,\alpha}}\leq c\left((T+t)^{-1-\frac{1}{p-2}}+\|w^\top\|^{p-1}_{C^{0,\alpha}}\right).
\end{equation}
Combining (\ref{eq27}), \eqref{taylor1}, and \eqref{taylor2}, we have
\begin{equation}\label{eq10}
  \|DE(u)\|_{C^{0,\alpha}}\leq c\left((T+t)^{-1-\frac{1}{p-2}}+\|w^\top\|^{p-1}_{C^{0,\alpha}}+\|w^\perp\|_{C^{2,\alpha}}\right).
\end{equation}
We define
\begin{equation*}
  E_0^\top(w):=\textrm{proj}_{\Lambda_0}E_0(w),\ \ E_0^\perp(w):=\textrm{proj}_{\Lambda_0^\perp} E_0(w).
\end{equation*}
Now the asserted bounds for $E_0^\top(w)$ follow from the bound (\ref{eq9}) on $E_0(w)$, the estimate (\ref{eq10}) and the continuity of the map $\textrm{proj}_{\Lambda_0}:C^{0,\alpha}(M,g_\infty)\rightarrow \Lambda_0$,
\begin{equation*}
  \begin{split}
  \|\textrm{proj}_{\Lambda_0}f\|_{C^{0,\alpha}( M,g_\infty)}\leq c\,\|f\|_{C^{0,\alpha}( M,g_\infty)}.
\end{split}
\end{equation*}
 Note that this is a spatial bound, so it does not include the $t$-H\"{o}lder norm, but the desired space-time norm bound follows easily from it in the same spirit of \cite[Lemma 16]{Chodosh}. The bound for $E_0^\top(w_1)-E_0^\top(w_2)$ follows similarly. This together with the bound (\ref{eq9}) on $E_0(w)$ and the estimate (\ref{eq10})
gives the estimates for $E_0^\perp(w)$.
\end{proof}

As mentioned at the beginning of this section, we will reduce the weighted Yamabe flow to two flows, one on $\Lambda_0$ and the other on $\Lambda_0^\perp$. The following proposition explains how to do this.

\begin{prop}\label{prop17}
  There exists $T_0>0$, $\epsilon_0>0$ and $c>0$, all depending on $(g_\infty,\phi_\infty)$ and $\hat{v}$, such that the following holds: Fix $T>T_0$. Then, for $\varphi(t)$ as in Lemma \ref{lem15} and $w\in C^{2,\alpha}(M\times [0,\infty))$, there are functions $E^\top(w)$ and $E^\perp(w)$ such that $u:=1+\varphi+w^\top+\Phi(\varphi+w^\top)+w^\perp$ is a solution to the weighted Yamabe flow if and only if
\begin{align}
  \frac{8}{n+m-2}(w^\top)'+D^2F_p(\varphi)w^\top=&E^\top(w),\label{eq11}\\
(w^\perp)'-\mathcal{L}_\infty w^\perp=&E^\perp(w).\label{eq12}
\end{align}
Here, as long as $\|w\|_{C^{2,\alpha}}\leq \epsilon_0$, the error terms $E^\top$ and $E^\perp$ satisfy
\begin{equation*}
  \begin{split}
  \|E^\top(w)\|_{C^{0,\alpha}}\leq & c\left((T+t)^{-\frac{p-1}{p-2}}+\|w^\top\|^{p-1}_{C^{0,\alpha}}+\|w^\perp\|_{C^{2,\alpha}}\right)\left((T+t)^{-\frac{1}{p-2}}+\|w\|_{C^{2,\alpha}}\right)\\
&+c(T+t)^{-\frac{p}{p-2}}+c(T+t)^{-\frac{p-1}{p-2}}\|w^\top\|_{C^{0,\alpha}}+c(T+t)^{-\frac{p-3}{p-2}}\|w^\top\|^2_{C^{0,\alpha}}\\
&+c\|w^\top\|^{p-1}_{C^{0,\alpha}}+c\left((T+t)^{-\frac{1}{p-2}}+\|w\|_{C^{2,\alpha}}\right)\|w^\perp\|_{C^{2,\alpha}},
\end{split}
\end{equation*}
\begin{equation*}
  \begin{split}
 & \|E^\top(w_1)-E^\top(w_2)\|_{C^{0,\alpha}}\\
&\qquad \leq c\left((T+t)^{-\frac{p-1}{p-2}}+\|w_1^\top\|^{p-1}_{C^{0,\alpha}}+\|w_2^\top\|^{p-1}_{C^{0,\alpha}}+\|w_1^\perp\|_{C^{2,\alpha}}+\|w_2^\perp\|_{C^{2,\alpha}}\right)\\
&\qquad\qquad \times\|w_1-w_2\|_{C^{2,\alpha}}\\
&\qquad\quad +c\left((T+t)^{-\frac{1}{p-2}}+\|w_1\|_{C^{2,\alpha}}+\|w_2\|_{C^{2,\alpha}}\right)(\|w_1^\top\|^{p-2}_{C^{0,\alpha}}+\|w_2^\top\|^{p-2}_{C^{0,\alpha}})\\
&\qquad\qquad \times \|w_1^\top-w_2^\top\|_{C^{0,\alpha}}\\
&\qquad \quad +c\left((T+t)^{-\frac{1}{p-2}}+\|w_1\|_{C^{2,\alpha}}+\|w_2\|_{C^{2,\alpha}}\right)\|w_1^\perp-w_2^\perp\|_{C^{2,\alpha}}\\
&\qquad\quad +c\left((T+t)^{-\frac{p-3}{p-2}}(\|w_1^\top\|_{C^{0,\alpha}}+\|w_2^\top\|_{C^{0,\alpha}})+\|w_1^\top\|^{p-2}_{C^{0,\alpha}}+\|w_2^\top\|^{p-2}_{C^{0,\alpha}}\right)\\
&\qquad \qquad \times \|w_1^\top-w_2^\top\|_{C^{0,\alpha}}\\
&\qquad\quad +c(T+t)^{-\frac{p-1}{p-2}}\|w_1^\top-w_2^\top\|_{C^{0,\alpha}},
\end{split}
\end{equation*}
\begin{equation*}
  \begin{split}
&  \|E^\perp(w)\|_{C^{0,\alpha}}\\
&\qquad \quad \leq c\left((T+t)^{-\frac{p-1}{p-2}}+\|w^\top\|^{p-1}_{C^{0,\alpha}}+\|w^\perp\|_{C^{2,\alpha}}\right)\left((T+t)^{-\frac{1}{p-2}}+\|w\|_{C^{2,\alpha}}\right)\\
&\qquad\qquad +c\left((T+t)^{-\frac{1}{p-2}}+\|w\|_{C^{2,\alpha}}\right)\|w^\perp\|_{C^{2,\alpha}}\\
&\qquad\qquad +c\left((T+t)^{-\frac{1}{p-2}}+\|w\|_{C^{2,\alpha}}\right)\left((T+t)^{-\frac{p-1}{p-2}}+\|w'\|_{C^{0,\alpha}}\right)
\end{split}
\end{equation*}
\begin{equation*}
  \begin{split}
 & \|E^\perp(w_1)-E^\perp(w_2)\|_{C^{0,\alpha}}\\
&\qquad\quad \leq c\left((T+t)^{-\frac{p-1}{p-2}}+\|w_1^{\top}\|^{p-1}_{C^{0,\alpha}}+\|w_2^\top\|^{p-1}_{C^{0,\alpha}}+\|w_1^\perp\|_{C^{2,\alpha}}+\|w_2^\perp\|_{C^{2,\alpha}}\right)\\
&\qquad\qquad \times\|w_1-w_2\|_{C^{2,\alpha}}\\
&\qquad\quad +c\left((T+t)^{-\frac{1}{p-2}}+\|w_1\|_{C^{2,\alpha}}+\|w_2\|_{C^{2,\alpha}}\right)(\|w_1^\top\|^{p-2}_{C^{0,\alpha}}+\|w_2^\top\|^{p-2}_{C^{0,\alpha}})\\
&\qquad\qquad \times \|w_1^\top-w_2^\top\|_{C^{0,\alpha}}\\
&\qquad\quad +c\left((T+t)^{-\frac{1}{p-2}}+\|w_1\|_{C^{2,\alpha}}+\|w_2\|_{C^{2,\alpha}}\right)\|w_1^\perp-w_2^\perp\|_{C^{2,\alpha}}\\
&\qquad\quad +c\left((T+t)^{-\frac{1}{p-2}}+\|w_1\|_{C^{2,\alpha}}+\|w_2\|_{C^{2,\alpha}}\right)\|w_1'-w_2'\|_{C^{0,\alpha}}\\
&\qquad\quad +c\left((T+t)^{-\frac{p-1}{p-2}}+\|w_1'\|_{C^{0,\alpha}}+\|w_2'\|_{C^{0,\alpha}}\right)\|w_1-w_2\|_{C^{2,\alpha}}.
\end{split}
\end{equation*}
Here we are using the parabolic H\"{o}lder norms on $(t,t+1)\times M$ as defined above; the bounds hold for each fixed $t\geq 0$, with the constants independent of $T$ and $t$.
\end{prop}
\begin{proof}
  Recall that $u$ is a solution to the weighted Yamabe flow if and only if
\begin{equation}\label{reduce1}
\begin{split}
 \frac{4}{n+m-2} \frac{\partial}{\partial t}u(t)=(r_{\phi(t)}^m-R_{\phi(t)}^m)u(t)
  \end{split}
\end{equation}
where $(g(t),\phi(t))=(u(t)^\frac{4}{m+n-2}g_\infty, \phi_\infty-\frac{2m}{m+n-2}\ln u(t))$. Regarding $DE(w)$ in \eqref{diffE} as an element of $L^2(M,e^{-\phi_\infty}dV_{g_\infty})$, we have
\begin{equation}\label{reduce2}
  \frac{1}{2}DE(w)=-\frac{4(n+m-1)}{n+m-2}\Delta_{\phi_\infty}w+R^m_{\phi_\infty}w-r^m_{\phi}w^{\frac{n+m+2}{n+m-2}}
\end{equation}
where as always in this paper, $E$ is defined on the unit volume conformal class. Combining \eqref{reduce1} and \eqref{reduce2}, we have
\begin{equation*}
  \frac{4}{n+m-2}\frac{\partial u}{\partial t}=-\frac{1}{2}DE(u)u^{-\frac{m+n+2}{m+n-2}}.
  \end{equation*}

We now project the weighted Yamabe flow equation onto $\Lambda_0$ and $\Lambda_0^\perp$, so $u$ solves the weighted Yamabe flow if and only if the following two equations are satisfied:
\begin{equation}\label{decompose}
  \begin{split}
  &\frac{8}{n+m-2}(\varphi+w^\top)'\\
&\qquad\qquad=-\textrm{proj}_{\Lambda_0}\left[DE\left(1+\varphi+w^\top+\Phi(\varphi+w^\top)+w^\perp\right)\right]-E_0^\top(w),\\
  &\frac{8}{n+m-2}\left(\Phi(\varphi+w^\top)+w^\perp\right)'\\
&\qquad\qquad=-\textrm{proj}_{\Lambda_0^\perp}\left[DE\left(1+\varphi+w^\top+\Phi(\varphi+w^\top)+w^\perp\right)\right]-E_0^\perp(w),
\end{split}
\end{equation}
where $E_0(w)$ is defined as in Lemma \ref{lemma16}.
Now, by the Taylor's theorem we claim that
\begin{equation}\label{proj1}
  \begin{split}
  &\textrm{proj}_{\Lambda_0}DE\left(1+\varphi+w^\top+\Phi(\varphi+w^\top)+w^\perp\right)\\
&\qquad\qquad\qquad =\textrm{proj}_{\Lambda_0}DE\left(1+\varphi+w^\top+\Phi(\varphi+w^\top)\right)+E_1^\top(w),
\end{split}
\end{equation}
with the bounds
\begin{equation}\label{bound1}
  \begin{split}
   \|E_1^\top(w)\|_{C^{0,\alpha}}\leq &c\left((T+t)^{-\frac{1}{p-2}}+\|w\|_{C^{2,\alpha}}\right)\|w^\perp\|_{C^{2,\alpha}},\\
\|E_1^\top(w_1)-E_1^\top(w_2)\|_{C^{0,\alpha}}\leq &c\left(\|w_1\|_{C^{2,\alpha}}+\|w_2\|_{C^{2,\alpha}}\right)\|w_1^\perp-w_2^\perp\|_{C^{2,\alpha}}.
\end{split}
\end{equation}
The claim follows from the integral form of the remainder in Taylor's theorem (see \cite[Proposition 17]{Chodosh} for more details) so we omit the proof here.
\medskip

Recall that $F(v):=E(\Psi(v))=E(1+v+\Phi(v))$, and using the Lyapunov-Schmidt reduction (Proposition \ref{prop7})
\begin{equation}\label{proj2}
  \textrm{proj}_{\Lambda_0}DE\left(1+\varphi+w^\top+\Phi(\varphi+w^\top)\right)=DF(\varphi+w^\top).
\end{equation}
 Furthermore, by analyticity (Proposition \ref{prop7}), $DF$ has a convergent power series representation around $0$ with lowest
 order term of order $p-1$. Thus, as long as $\varphi+w^\top$ is small enough, we may write
\begin{equation}\label{proj3}
  DF(\varphi+w^\top)=DF(\varphi)+D^2F(\varphi)(w^\top)+E_2^\top(w^\top),
\end{equation}
where
\begin{equation}\label{bound2}
  \begin{split}
  &\|E_2^\top(w^\top)\|_{C^{0,\alpha}}\leq c\left((T+t)^{-\frac{p-3}{p-2}}+\|w^\top\|^{p-3}_{C^{0,\alpha}}\right)\|w^\top\|^2_{C^{0,\alpha}},\\
  &\|E_2^\top(w_1^\top)-E_2^\top(w_2^\top)\|_{C^{0,\alpha}}\\
&\qquad\qquad \leq c\left((T+t)^{-\frac{p-3}{p-2}}(\|w_1^\top\|_{C^{0,\alpha}}+\|w_2^\top\|_{C^{0,\alpha}})+\|w_1^\top\|^{p-2}_{C^{0,\alpha}}+\|w_2^\top\|^{p-2}_{C^{0,\alpha}}\right)\\
&\qquad\qquad \qquad \times\|w_1^\top-w_2^\top\|_{C^{0,\alpha}}.
\end{split}
\end{equation}
\medskip

By the results we have obtained so far, the $\Lambda_0$-component of the weighted Yamabe flow may be written as
\begin{equation*}
  \begin{split}
  \frac{8}{n+m-2}(\varphi+w^\top)'
=-DF(\varphi)-D^2F(\varphi)(w^\top)-E_2^\top(w^\top)-E_1^\top(w)-E_0^\top(w)
\end{split}
\end{equation*}
where we have used (\ref{proj1}), (\ref{proj2}) and (\ref{proj3}). Now, expanding $F$ in a power series, $F=F(0)+\sum_{j=p}^\infty F_j$, we may write the above expression as
\begin{equation*}
  \begin{split}
  &\frac{8}{n+m-2}(\varphi+w^\top)'\\
&\qquad\qquad=-DF_p(\varphi)-D^2F_p(\varphi)(w^\top)+\underbrace{E_3^\top(w)-E_2^\top(w^\top)-E_1^\top(w)-E_0^\top(w)}_{:=E^\top(w)}
\end{split}
\end{equation*}
where
\begin{equation*}
  E_3^\top(w)=\sum_{j\geq p+1}(DF_j(\varphi)+D^2F_j(\varphi)w^\top).
\end{equation*}
On the other hand, by Lemma \ref{lem15}, we have that $\frac{8}{n+m-2}\varphi'=-DF_p(\varphi)$. Therefore, $w^\top$ must satisfy the equation
\begin{equation*}
  \frac{8}{n+m-2}(w^\top)'+D^2F_p(\varphi)w^\top=E^\top(w).
\end{equation*}
By analyticity, $E_3^\top(w)$ converges in $C^{0,\alpha}$ for $\|\varphi\|_{C^{2,\alpha}}$ and $\|w\|_{C^{2,\alpha}}$ small enough. Because each term in the sum is a homogeneous polynomial, we get the following error bound by using the formula for $\varphi$:
\begin{equation}\label{bound3}
\begin{split}
    \|E_3^\top(w)\|_{C^{0,\alpha}}&\leq c\left((T+t)^{-\frac{p}{p-2}}+(T+t)^{-\frac{p-1}{p-2}}\|w^\top\|_{C^{0,\alpha}}\right),\\
\|E_3^\top(w_1)-E_3^\top(w_2)\|_{C^{0,\alpha}}&\leq c(T+t)^{-\frac{p-1}{p-2}}\|w_1^\top-w_2^\top\|_{C^{0,\alpha}}.
\end{split}
\end{equation}
Combining (\ref{bound1}), (\ref{bound2}) and (\ref{bound3}), we can see that $E^\top(w)$ satisfies the asserted bounds.
By the similar argument, we can prove the result for the $\Lambda_0^{\perp}$-portion of the weighted Yamabe flow.
\end{proof}

\bigskip

\subsection{Solving the kernel-projected flow with polynomial decay estimates}\label{sec3.2}

In this subsection we solve the kernel-projected flow (\ref{eq11}). First, from the definition of $\varphi$ in (\ref{eq6}) and the fact that $D^2F_p$ is homogeneous of degree $p-2$,
\begin{equation*}
  D^2F_p(\varphi)=(T+t)^{-1}\underbrace{\left(\frac{8}{(n+m-2)p(p-2)F_p(\hat{v})}\right)D^2F_p(\hat{v})}_{:=\mathcal{D}}.
\end{equation*}
Let $\mu_1,\cdots,\mu_k$ be the eigenvalues of $\mathcal{D}$ and $e_i$ the corresponding orthonormal basis in which $\mathcal{D}$ is diagonalized. Then the kernel-projected flow is equivalent to the following system of ODEs for $v_i:=w^\top\cdot e_i$,
\begin{equation}\label{eq14}
  \frac{8}{n+m-2}v_i'+\frac{\mu_i}{T+t}v_i=E_i^\top:=E^\top\cdot e_i,\quad i=1,\cdots,k.
\end{equation}
Fix for the rest of this subsection a number $\gamma$ with $\gamma\notin\left\{\frac{n+m-2}{8}\mu_1,\cdots,\frac{n+m-2}{8}\mu_k\right\}$.
Define the following weighted norms:
\begin{equation*}
  \|u\|_{C^{0,\alpha}_\gamma}:=\sup_{t>0}\left[(T+t)^\gamma\|u(t)\|_{C^{0,\alpha}}\right] \ \textrm{ and } \ \|u\|_{C^{0,\alpha}_{1,\gamma}}:=\|u\|_{C^{0,\alpha}_\gamma}+\|u'\|_{C^{0,\alpha}_{1+\gamma}}.
\end{equation*}
We recall that these H\"{o}lder norms are space-time norms on the interval $(t,t+1)\times  M$, as defined in (\ref{eq8}).

Given $\gamma$ as above, we define $\Pi_0=\Pi_0(\gamma)$ by
\begin{equation}
  \Pi_0:=\textrm{span}\left\{v\in\Lambda_0:\mathcal{D}v=\mu v,\textrm{ and }\mu>\frac{8}{n+m-2}\gamma\right\}.
\end{equation}
Moreover, let $\textrm{proj}_{\Pi_0}:\Lambda_0\rightarrow \Pi_0$ be the corresponding linear projector. The next lemma concerns the system (\ref{eq14}).

\begin{lem}\label{lemma18}
  For any $T>0$ such that $\|E^\top\|_{C^{0,\alpha}_{1+\gamma}}<\infty$, there is a unique $u$ with $u(t)\in\Lambda_0$, $t\in[0,\infty)$, satisfying $\|u\|_{C^0_\gamma}<\infty$, $\textrm{\emph{proj}}_{\Pi_0}(u(0))=0$, and such that $v_i:=u\cdot e_i$ solves the system \eqref{eq14}. Furthermore, we have the bound
\begin{equation*}
  \|u\|_{C^{0,\alpha}_{1,\gamma}}\leq C\|E^\top\|_{C^{0,\alpha}_{1+\gamma}}.
\end{equation*}
Here the constant $C$ does not depend on $T$.
\end{lem}
\begin{proof}
Letting
\begin{equation*}
  w_j:=(T+t)^{\frac{n+m-2}{8}\mu_j}v_j,
\end{equation*}
the system (\ref{eq14}) is equivalent to
\begin{equation*}
  w_j'=\frac{n+m-2}{8}(T+t)^{\frac{n+m-2}{8}\mu_j}E_j^\top,\quad j=1,\cdots,k.
\end{equation*}
Then, we claim that we may solve the $j$-th ODE as
\begin{equation*}
  w_j(t)=\left\{
           \begin{array}{ll}
             \displaystyle\alpha_j-\frac{n+m-2}{8}\int_t^\infty(T+\tau)^{\frac{n+m-2}{8}\mu_j}E_j^\top(\tau)d\tau, & \hbox{if $\gamma>\displaystyle\frac{n+m-2}{8}\mu_j$;} \\
             \displaystyle\alpha_j+\frac{n+m-2}{8}\int_0^t(T+\tau)^{\frac{n+m-2}{8}\mu_j}E_j^\top(\tau)d\tau, & \hbox{if $\gamma<\displaystyle\frac{n+m-2}{8}\mu_j$.}
           \end{array}
         \right.
\end{equation*}

First suppose that $j$ is such that $\gamma>\frac{n+m-2}{8}\mu_j$. Then the claim would imply that
\begin{equation*}
\begin{split}
   v_j(t)=(T+t)^{-\frac{n+m-2}{8}\mu_j}\alpha_j-\frac{n+m-2}{8}(T+t)^{-\frac{n+m-2}{8}\mu_j}\int_t^\infty(T+\tau)^{\frac{n+m-2}{8}\mu_j}E_j^\top(\tau)d\tau.
\end{split}
\end{equation*}
To prove the claim, we check that the integral converges under our assumption on $E^\top$:
\begin{equation*}
  \begin{split}
  &\left|(T+t)^{-\frac{n+m-2}{8}\mu_j}\int_t^\infty(T+\tau)^{\frac{n+m-2}{8}\mu_j}E_j^\top(\tau)d \tau\right|\\
  &\qquad\leq(T+t)^{-\frac{n+m-2}{8}\mu_j}\|E_j\|^\top_{C^0_{1+\gamma}}\int_t^\infty(T+\tau)^{\frac{n+m-2}{8}\mu_j-\gamma-1}d\tau\\
  &\qquad=\left(\gamma-\frac{n+m-2}{8}\right)^{-1}(T+t)^{-\frac{n+m-2}{8}\mu_j}\|E_j^\top\|_{C^0_{1+\gamma}}(T+t)^{\frac{n+m-2}{8}\mu_j-\gamma}\\
  &\qquad = C_j(T+t)^{-\gamma}\|E_j^\top\|_{C^0_{1+\gamma}}.
\end{split}
\end{equation*}
The previous estimate also shows that, since by assumption $\gamma>\frac{n+m-2}{8}\mu_j$, to have $\|u\|_{C^{0}_\gamma}<\infty$, it must hold that $\alpha_j=0$.

On the other hand, if $\gamma<\frac{n+m-2}{8}\mu_j$, by requiring $\textrm{proj}_{\Pi_0}u(0)=0$, we see that $\alpha_j=0$. As a result, the bounds for $\|v_j\|_{C^0_\gamma}$ follow from a similar calculation as before. Combining these two cases proves existence, uniqueness and the $\|u\|_{C^0_\gamma}$ bound.

 It thus remains to prove the inequality $\|u\|_{C^{0,\alpha}_{1,\gamma}}\leq C\|E^\top\|_{C^{0,\alpha}_{1+\gamma}}$.
By finite dimensionality, the (spatial) $C^{0,\alpha}(M)$-H\"{o}lder norms of each basis element in $\Lambda_0$ are uniformly bounded. Thus, it remains to show that the desired inequality holds for the H\"{o}lder norms in the time direction, along with the same thing for $u'(t)$.

Suppose that $j$ is such that $\gamma>\frac{n+m-2}{8}\mu_j$. Then, we have seen above that
$$
v_j(t)=-\frac{n+m-2}{8}(T+t)^{-\frac{n+m-2}{8}\mu_j}\int_t^\infty(T+\tau)^{\frac{n+m-2}{8}\mu_j}E_j^\top(\tau)d\tau,$$
which gives
\begin{equation*}
\begin{split}
   v_j'(t)=&\frac{(n+m-2)^2}{64}\mu_j(T+t)^{-\frac{n+m-2}{8}\mu_j-1}\int_t^\infty(T+\tau)^{\frac{n+m-2}{8}\mu_j}E_j^\top(\tau)d \tau\\
   &+\frac{n+m-2}{8}E_j^\top(t).
\end{split}
\end{equation*}
Thus
\begin{equation*}
  \begin{split}
  \|v_j'\|_{C^{0,\alpha}}&\leq C\left\|(T+t)^{-\frac{n+m-2}{8}\mu_j-1}\int_t^\infty(T+\tau)^{\frac{n+m-2}{8}\mu_j}E_j^\top(\tau)d\tau\right\|_{C^1}+C\|E_j^\top(\tau)\|_{C^{0,\alpha}}\\
  &\leq C\left\|(T+t)^{-\frac{n+m-2}{8}\mu_j-2}\int_t^\infty(T+\tau)^{\frac{n+m-2}{8}\mu_j}E_j^\top(\tau)d\tau\right\|_{C^0}+C\|E_j^\top(\tau)\|_{C^{0,\alpha}}\\
  &\leq C(T+t)^{-1-\gamma}\|E_j^\top\|_{C^{0,\alpha}_{1+\gamma}}.
\end{split}
\end{equation*}
On the other hand, the case of $\gamma<\frac{n+m-2}{8}\mu_j$ can be easily be obtained through a similar argument. Combining these calculations, we obtain a H\"{o}lder estimate for $v_j$. From this the claimed inequality follows.
\end{proof}

\subsection{Solving the kernel-orthogonal projected flow}\label{sec3.3}
In this subsection, we solve the kernel-orthogonal projected flow, which is the remaining part of the weighted Yamabe flow.
Define the weighted norms
\begin{equation*}
  \|u\|_{L^2_q}=\sup_{t\in[0,\infty)}\left[(T+t)^q\|u(t)\|_{L^2( M)}\right],
\end{equation*}
where the $L^2$ norm is the spatial norm of $u(t)$ on $ M$, taken with respect to $e^{-\phi_\infty}dV_{g_\infty}$, and
\begin{equation*}
  \|u\|_{C^{2,\alpha}_q}=\sup_{t\geq 0} \left[(T+t)^q\|u(t)\|_{C^{2,\alpha}}\right],
\end{equation*}
where, as usual, the H\"{o}lder norms are the space-time norms defined in (\ref{eq8}). Also, let
\begin{equation*}
  \begin{split}
  &\Lambda_\downarrow:=\overline{\textrm{span}\{\varphi\in C^\infty( M):\mathcal{L}_\infty \varphi+\delta \varphi=0,\delta>0\}}^{L^2},\\
&\Lambda_\uparrow :=\textrm{span}\{\varphi\in C^\infty( M):\mathcal{L}_\infty \varphi+\delta \varphi=0,\delta<0\}.
\end{split}
\end{equation*}

From the spectral theory, $L^2( M,e^{-\phi_\infty} g_\infty)=\Lambda_\uparrow\oplus\Lambda_0\oplus \Lambda_\downarrow$ and $\Lambda_\uparrow$ and $\Lambda_0$ are finite dimensional. Write the nonnegative integers as an ordered union $\mathbb{N}=K_\uparrow\cup K_0\cup K_\downarrow$, where the ordering of the indices comes from an ordering of the eigenfunctions of the $\mathcal{L}_\infty$ and the partitioning of $\mathbb{N}$ corresponds to which of $\Lambda_\downarrow$, $\Lambda_0$, or $\Lambda_\uparrow$ the $k$-th eigenfunction of $\mathcal{L}_\infty$ lies in.

\begin{lem}\label{lemma19}
  For any $T>0$ and $q<\infty$ such that $\|E^\perp\|_{L^2_q}<\infty$, there is a unique $u(t)$ with $u(t)\in\Lambda_0^\perp$, $t\in[0,\infty)$, satisfying $\|u\|_{L^2_q}<\infty$, $\textrm{\emph{proj}}_{\Lambda_\downarrow}(u(0))=0$ and
\begin{equation}\label{eq15}
  u'=\mathcal{L}_\infty u+E^\perp.
\end{equation}
Furthermore, $\|u\|_{L^2_q}\leq C\|E^\perp\|_{L^2_q}$ and $\|u\|_{C^{2,\alpha}_q}\leq C\|E^\perp\|_{C^{0,\alpha}_q}$. Here the constant $C$ does not depent on $T$.
\end{lem}
\begin{proof}
Let $\varphi_i$ be an eigenfunction of $\mathcal{L}_\infty$ with eigenvalue $-\delta_i$ which is orthogonal to the kernel $\Lambda_0$. The equation (\ref{eq15}) reduces to the system
\begin{equation}\label{eq17}
  u_i'+\delta_iu_i=E_i^\perp,
\end{equation}
where $u_i=\langle u,\varphi_i\rangle$ and $E_i^\perp=\langle E^\perp,\varphi_i\rangle$. This is equivalent to
\begin{equation}\label{eq18}
  (e^{\delta_it}u_i)'=e^{\delta_i t}E_i^\perp.
\end{equation}
Thus, we may represent the components of the solution as
\begin{equation*}
\begin{split}
u_i(t)=\begin{dcases}
  &\beta_ie^{-\delta_it}+e^{-\delta_it}\int_0^te^{\delta_i \tau}E_i^\perp(\tau)d\tau~~\mbox{ for }i\in K_\downarrow,\\
&\beta_ie^{-\delta_it}-e^{-\delta_it}\int_t^\infty e^{\delta_i \tau}E_i^\perp(\tau)d\tau~~\mbox{ for }i\in K_{\uparrow}.
\end{dcases}
\end{split}
\end{equation*}
In particular, we have
\begin{equation*}
\begin{split}
  u(t)&=\sum_{j\in K_\downarrow}\left(\beta_je^{-\delta_jt}+e^{-\delta_jt}\int_0^te^{\delta_j \tau}E_j^\perp(\tau)d\tau\right)\varphi_j\\
&\qquad+\sum_{j\in K_\uparrow}\left(\beta_je^{-\delta_jt}-e^{-\delta_jt}\int_t^\infty e^{\delta_j \tau}E_j^\perp(\tau)d\tau\right)\varphi_j.
\end{split}
\end{equation*}
This sum is in an $L^2$ sense (but then elliptic regularity guarantees that the sum converges uniformly on compact time intervals). We note that for $i\in K_\uparrow$, if $\|u\|_{L^2_q}<\infty$, then necessarily $\beta_i=0$. Furthermore, by requiring that $\textrm{proj}_{\Lambda_\downarrow}u(0)=0$, we also have $\beta_i=0$ for $i\in K_\downarrow$.

The $L^2$-bound for the first term in $u$  can be estimated as:
\begin{equation*}
  \begin{split}
  \Bigg\|\sum_{j\in K_\downarrow}u_j(t)\varphi_j\Bigg\|_{L^2}^2&\leq \sum_{j\in K_\downarrow}\left(\int_0^t e^{\delta_j( \tau-t)}E_j^\perp(\tau)d\tau\right)^2\\
  &\leq \sum_{j\in K_\downarrow}\left(\int_0^t e^{\delta_{\min}( \tau-t)}E_j^\perp(\tau)d\tau\right)^2\\
  &\leq \left\|\int_0^t e^{\delta_{min}(\tau-t)}E^\perp(\tau)d\tau\right\|_{L^2}^2
\end{split}
\end{equation*}
where $\delta_{min}=\min_{j\in K_\downarrow}\delta_j$ and the inequality follows from the Parseval identity. Taking square roots and using the decay assumption on $E^\perp$ gives
\begin{equation*}
\begin{split}
   \Bigg\|\sum_{j\in K_\downarrow}u_j(t)\varphi_j\Bigg\|_{L^2}&\leq \left\|\int_0^t e^{\delta_{min}(\tau-t)}E^\perp(\tau)d\tau\right\|_{L^2}\\
   & \leq \int_0^t e^{\delta_{min}(\tau-t)}\left\|E^\perp(\tau)\right\|_{L^2}d\tau\\
   & \leq \left\|E^\perp\right\|_{q}\int_0^t e^{\delta_{min}(\tau-t)}(T+\tau)^{-q}d\tau\\
   & \leq C\|E^\perp\|_{L^2_q}(T+t)^{-q}.
   \end{split}
\end{equation*}
A similar argument holds for the $K_\uparrow$ terms. From this, the asserted bounds for $\|u\|_{L^2_q}$ follow readily.

We now consider the $C^{2,\alpha}_q$ bounds for $u$. Following the argument in \cite[Lemma 19]{Chodosh}, by interior parabolic Schauder estimates \cite[Theorem 4.9]{Lieberman} and Arzel\`{a}-Ascoli theorem, we have that for $t\geq 1$,
\begin{equation*}
  \begin{split}
  &\|u(t)\|_{C^{2,\alpha}}\leq C\left(\sup_{(s,x)\in(t-1,t+1)\times M}\|u(s,x)\|_{L^2(M)}+\|E^\perp\|_{C^{0,\alpha}((t-1,t+1)\times  M)}\right)\\
&\hspace{6cm}  +C\epsilon \|u(t)\|_{C^{0,\alpha}((t-1,t+1)\times  M)}.
\end{split}
\end{equation*}
Multiplying it by $(T+t)^q$ and taking the supremum over $t\geq1$ yields
\begin{equation}\label{4.3.1}
\begin{split}
 & \sup_{t\geq 1}\left[(T+t)^q\|u(t)\|_{C^{2,\alpha}}\right]\leq C\|E^\perp\|_{C^{0,\alpha}_q}+C\epsilon\|u\|_{C^{0,\alpha}_q}
\end{split}
\end{equation}
where we have used  $\|u\|_{L^2_q}\leq C\|E^\perp\|_{L^2_q}$, which was proved earlier.
To finish the proof, it remains to extend the supremum up to $t=0$. The global Schauder estimates \cite[Theorem 4.28]{Lieberman} shows that
\begin{equation}\label{4.3.2}
  \begin{split}
  \|u(t)\|_{C^{2,\alpha}((0,1)\times M)}&\leq C\Big(\sup_{s\in(0,1)}\|u(s,x)\|_{L^2( M)}+\epsilon\|u\|_{C^{0,\alpha}((0,1)\times M)}\\
&\qquad\qquad +\|E^\perp\|_{C^{0,\alpha}((0,1)\times  M)}+\|u(0)\|_{C^{2,\alpha}(M)}\Big).
\end{split}
\end{equation}
Except for the last term  $\|u(0)\|_{C^{2,\alpha}( M)}$ on the right-hand side of the above expression, the rest of the terms can be bounded in a manner similar to the argument used above. Note that
\begin{equation*}
  u(0)=-\sum_{j\in K_\uparrow}\left(\int_{0}^\infty e^{\delta_j\tau}E_j^\perp(\tau)d\tau\right)\varphi_j.
\end{equation*}
The space $\Lambda_\uparrow$ is finite-dimensional, so there must be a uniform constant $C>0$ such that $\|\varphi_j\|_{C^{2,\alpha}( M)}\leq C\|\varphi_j\|_{L^2( M)}$ for all $j\in K_\uparrow$. Using this we have that
\begin{equation}\label{4.3.3}
  \begin{split}
  \|u(0)\|_{C^{2,\alpha}( M)}^2& \leq  C\sum_{j\in K_\uparrow}\left(\int_0^\infty e^{\delta_j\tau}E_j^\perp(\tau)d\tau\right)^2\|\varphi_j\|^2_{C^{2,\alpha}(M)}\\
  & \leq  C\sum_{j\in K_\uparrow}\left(\int_0^\infty e^{\delta_j\tau}E_j^\perp(\tau)d\tau\right)^2\|\varphi_j\|^2_{L^2( M)}=C\|u(0)\|^2_{L^2( M)}.
\end{split}
\end{equation}
Combining (\ref{4.3.1}), (\ref{4.3.2}), and (\ref{4.3.3}), we obtain that
\begin{equation*}
  \sup_{t\geq 0}\left[(T+t)^q\|u(t)\|_{C^{2,\alpha}}\right]\leq C\|E^\perp\|_{C^{0,\alpha}_q}+C\epsilon\|u\|_{C^{0,\alpha}_q}.
\end{equation*}
By choosing $\epsilon$ small, we get the desired H\"{o}lder bounds.
\end{proof}

\subsection{Construction of a slowly converging flow}\label{sec3.4}
In this subsection we will combine the results from the previous two subsections to finally construct a slowly converging flow. That is, we prove Theorem \ref{main2}.

To proceed further, we define the norm
\begin{equation*}
  \|f\|_{\gamma}^*:=\|\textrm{proj}_{\Lambda_0}f\|_{C^{0,\alpha}_{1,\gamma}}+\|\textrm{proj}_{\Lambda_0^\perp}f\|_{C^{2,\alpha}_{1+\gamma}}.
\end{equation*}
Recall that
\begin{equation*}
  \begin{split}
  &\|u\|_{C^{0,\alpha}_{1,\gamma}}=\sup_{t\geq 0}\left[(T+t)^\gamma\|u(t)\|_{C^{0,\alpha}}\right]+\sup_{t\geq0}\left[(T+t)^{1+\gamma}\|u'(t)\|_{C^{0,\alpha}}\right],\\
&\|u\|_{C^{2,\alpha}_{1+\gamma}}=\sup_{t\geq 0}\left[(T+t)^{1+\gamma}\|u(t)\|_{C^{2,\alpha}}\right],
\end{split}
\end{equation*}
where the H\"{o}lder norms are the space-time H\"{o}lder norms defined in (\ref{eq8}). For $\gamma$ to be specified below,  the Banach space $X$ is defined as
\begin{equation}\label{defX}
  X:=\{w:\|w\|_{\gamma}^*<\infty\}.
\end{equation}

\begin{prop}\label{prop4.7}
  Assume that $(g_\infty,\phi_\infty)$ satisfies $AS_p$. We may thus fix a point where $F_p|_{\mathbb{S}^{k-1}}$ achieves a positive maximum and denote it by $\hat{v}$. Define
\begin{equation*}
  \varphi(t)=(T+t)^{-\frac{1}{p-2}}\left(\frac{8}{(n+m-2)p(p-2)F_p(\hat{v})}\right)^{\frac{1}{p-2}}\hat{v}
\end{equation*}
as in Lemma \ref{lem15}. Then, there exists $C>0$, $T>0$, $\frac{1}{p-2}<\gamma<\frac{2}{p-2}$ and $u(t)\in C^\infty(M\times (0,\infty))$ such that
  $u(t)>0$ for all $t>0$, $(g(t),\phi(t)):=(u(t)^\frac{4}{n+m-2}g_\infty,\phi_\infty-\frac{2m}{n+m-2}\ln u(t))$ is a solution to the weighted Yamabe flow and
\begin{equation*}
  \|w^\top(t)+\Phi\left(\varphi(t)+w^\top(t)\right)+w^\perp(t)\|_\gamma^*=\|u(t)-\varphi(t)-1\|_\gamma^*\leq C.
\end{equation*}
\end{prop}
\begin{proof}
  We fix $\frac{1}{p-2}<\gamma<\frac{2}{p-2}$ so that $\gamma\notin\left\{\frac{n+m-2}{8}\mu_1,\cdots,\frac{n+m-2}{8}\mu_k\right\}$. By Proposition \ref{prop17}, the weighted Yamabe flow can be reduced to two flows, i.e. kernel projected flow and kernel-orthogonal projected flow, so it is enough to solve
\begin{equation*}
  \frac{8}{n+m-2}(w^\top)'+D^2F_p(\varphi)w^\top=E^\top(w),\quad (w^\perp)'-\mathcal{L}_\infty w^\perp=E^\perp (w),
\end{equation*}
for $w(t)$ with $\|w\|^*_\gamma<C$. To do so, we will use the contraction mapping method. We define a map
\begin{equation*}
  S:\{w\in X:\|w\|_\gamma^*\leq1\}\rightarrow X
\end{equation*}
where $X$ is the Banach space defined in (\ref{defX}), by defining $u:=\textrm{proj}_{\Lambda_0}S(w)$ and $v:=\textrm{proj}_{\Lambda_0^\perp}S(w)$ to be the solution of the kernel projected flow and the kernel-orthogonal projected flow with the initial values $u(0)=\textrm{proj}_{\Pi_0^\perp}w^\perp(0)$ and $v(0)=\textrm{proj}_{\Lambda_\uparrow}w^\perp(0)$ respectively, i.e.
\begin{equation*}
  \frac{8}{n+m-2}u'+D^2F_p(\varphi)u=E^\top(w)~~\mbox{ and }~~ v'-\mathcal{L}_\infty v=E^\perp (w).
\end{equation*}
Thus, we have defined the map $S(w)$ by its orthogonal projections onto $\Lambda_0$ and $\Lambda_0^\perp$. These solutions exist, in the right function spaces, by combining the bounds for the error terms in Proposition \ref{prop17} with Lemmas \ref{lemma18} and \ref{lemma19}. Furthermore, we have the explicit bound
\begin{equation*}
  \begin{split}
  \|\textrm{proj}_{\Lambda_0}S(w)\|_{C^{0,\alpha}_{1,\gamma}}&\leq c\|E^{\top}(w)\|_{C^{0,\alpha}_{1+\gamma}}\\
  & \leq c \sup_{t\geq 0} (T+t)^{1+\gamma}\left((T+t)^{-1-\frac{1}{p-2}}+\|w^{\top}\|^{p-1}_{C^{0,\alpha}}+\|w^\perp\|_{C^{2,\alpha}}\right)\\
  & \qquad \times \left((T+t)^{-\frac{1}{p-2}}+\|w\|_{C^{2,\alpha}}\right)\\
  & \qquad +c\sup_{t\geq 0} \left((T+t)^{\gamma-\frac{2}{p-2}}+(T+t)^{\gamma-\frac{2}{p-2}}\|w^\top\|_{C^{2,\alpha}}\right)\\
  & \qquad +c\sup_{t\geq 0} \left((T+t)^{\gamma+\frac{1}{p-2}}\|w^\top\|_{C^{2,\alpha}}+(T+t)^{\gamma+1}\|w^\top\|^{p-1}_{C^{2,\alpha}}\right)\\
  & \qquad +c\sup_{t\geq 0} (T+t)^{\gamma+1}\|w^\perp\|^{2}_{C^{2,\alpha}}\\
  & \leq c\left(T^{\gamma-\frac{2}{p-2}}+\left(T^{-\frac{1}{p-2}}+T^{(p-2)\left(\frac{1}{p-2}-\gamma\right)}\right)\|w\|^*_\gamma\right).
\end{split}
\end{equation*}

By the same argument, using Proposition \ref{prop17} and Lemma \ref{lemma19}, we obtain the similar bound for the kernel-orthogonal projected part:
\begin{equation*}
  \begin{split}
  \|\textrm{proj}_{\Lambda_0^\perp}S(w)\|_{C^{2,\alpha}_{1+\gamma}}\leq  c\left(T^{\gamma-\frac{2}{p-2}}+\left(T^{-\frac{1}{p-2}}+T^{(p-2)\left(\frac{1}{p-2}-\gamma\right)}\right)\|w\|^*_\gamma\right).
\end{split}
\end{equation*}
Therefore, we have
\begin{equation*}
\begin{split}
  \|S(w)\|_\gamma^*&=\|\textrm{proj}_{\Lambda_0}S(w)\|_{C^{0,\alpha}_{1,\gamma}}+\|\textrm{proj}_{\Lambda_0^\perp}S(w)\|_{C^{2,\alpha}_{1+\gamma}}\\
&\leq c\left\{T^{\gamma-\frac{2}{p-2}}+\left(T^{-\frac{1}{p-2}}+T^{(p-2)\left(\frac{1}{p-2}-\gamma\right)}\right)\|w\|^*_\gamma\right\}.
\end{split}
\end{equation*}
Thus, because $\gamma\in\left(\frac{1}{p-2},\frac{2}{p-2}\right)$, by choosing $T$ large enough we can ensure that $S$ maps $\{w:\|w\|_\gamma^*\leq 1\}\subset X$ into itself. And one can also show that $S$ is a contraction map by enlarging $T$ if necessary.
\end{proof}

Now we are ready to prove Theorem \ref{main2}.
\begin{proof}[Proof of Theorem \ref{main2}]
  From Proposition \ref{prop17}, we have constructed $\varphi(t)$ and $u(t)$ so that
\begin{equation*}
  \varphi(t)=(T+t)^{-\frac{1}{p-2}}\left(\frac{8}{(n+m-2)p(p-2)F_p(\hat{v})}\right)^{\frac{1}{p-2}}\hat{v},
\end{equation*}
 $(u(t)^\frac{4}{n+m-2}g_\infty,\phi_\infty-\frac{2m}{n+m-2}\ln u(t))$ is a solution to the weighted Yamabe flow, and
\begin{equation*}
  u(t)=1+\varphi(t)+\tilde{w}(t):=1+\varphi(t)+w^\top(t)+\Phi\left(\varphi(t)+w^\top(t)\right)+w^\perp(t),
\end{equation*}
where $\tilde{w}(t)$ satisfies $\|\tilde{w}\|_{C^0}\leq C(1+t)^{-\gamma}$ for some $C>0$ and all $t\geq 0$. We have arranged that $\gamma>1/(p-2)$, which implies that $\varphi(t)$ is decaying slower than $\tilde{w}(t)$. Thus
\begin{equation*}
  \|u(t)-1\|_{C^0}\geq C(1+t)^{-\frac{1}{p-2}}
\end{equation*}
as $t\rightarrow \infty$. From this, the assertion follows.
\end{proof}

\section{Examples satisfying $AS_p$}\label{section5}

In this section, we provide a family of metrics which satisfy $AS_3$. This allows us, via Theorem \ref{main2}, to conclude the existence of slowly converging weighted Yamabe flow.

We denote the $2n_2$-dimensional complex projective space equipped with the Fubini--Study metric by $(\P^{n_2}, g_{FS})$. (We normalize the Fubini--Study metric so the map $\S^{2n+1}\to \P^n$ from the unit sphere is a submersion.) From \cite[Section 5.1]{Chodosh}, we know that $R_{g_{FS}}=4n_2(n_2+1)$ and $\lambda_1(g_{FS})=4(n_2+1)$.

Suppose that $(M^{n_1}, g_M, e^{-\phi}dV_g, m)$ is a smooth metric measure space with constant weighted scalar curvature $R^m_{\phi}=\lambda_1(g_{FS})(n_1+n_2+m-1)$. We consider the product smooth metric measure space $(M^{n_1}\times \P^{n_2}, g=g_M\oplus g_{FS}, e^{-\phi}dV_{g_M\oplus g_{FS}}, m)$. Therefore, the weighted scalar curvature on it is given by
\begin{equation}
    R^m_{g,\phi}=R^m_{g_M,\phi}+R_{g_{FS}}=\lambda_1(n+m-1),
\end{equation}
where $n=n_1+2n_2$ is the dimension of the product space. So $\Lambda_0$ consists of eigenfunctions of $\Delta_{g,\phi}$ with eigenvalue $R^m_{g,\phi}/(n+m-1)=\lambda_1$. From this, we see that $(M^{n_1}\times \P^{n_2}, g=g_M\oplus g_{FS}, e^{-\phi}dV_{g_M\oplus g_{FS}}, m)$ is degenerate.

Moreover, let $v$ be an eigenfunction corresponding to the eigenvalue $\lambda_1$, i.e.
\begin{equation}\label{5.1}
-\Delta_{g_{FS}} v=\lambda_1v.
\end{equation}
It follows from the Courant nodal domain theorem that $v$ does not change sign. By replacing $v$ with $-v$ if necessary, we may assume that $v>0$.

Since $\phi\in C^\infty(M)$, the function $1\otimes v$ on $M^{n_1}\times \P^{n_2}$ will be an eigenfunction of $\Delta_{g,\phi}$ with eigenvalue $\lambda_1$, i.e.
\begin{displaymath}
(n+m-1)\Delta_{g,\phi}(1\otimes v)+R_{g,\phi}^m (1\otimes v)=(n+m-1)\Delta_{g_{FS}}v+(n+m-1)\lambda_1v=0.
\end{displaymath}

On the other hand, it follows from \cite[Section 5.1]{Chodosh} that
\begin{equation}
    \int_{\P^{n_2}}v^3 dV_{g_{FS}}\neq 0.
\end{equation}
Therefore, we have
\begin{equation*}
\begin{split}
F_3(1\otimes v)&=-2\left(\frac{n+m+2}{n+m-2}\right)\left(\frac{4}{n+m-2}\right)R_{g,\phi}^m\int_{M^{n_1}\times \P^{n_2}} v^3 e^{-\phi}dV_{g_M\oplus g_{FS}}\\
&=\frac{-8(n+m+2)(n+m-1)}{(n+m-2)^2}\lambda_1\left(\int_{M^{n_1}}e^{-\phi}dV_{g_M}\right)\left(\int_{\P^{n_2}}v^3dV_{g_FS}\right)\neq 0
\end{split}
\end{equation*}
since $v>0$ and $\lambda_1>0$.
Therefore, $(M^{n_1}\times \P^{n_2}, g=g_M\oplus g_{FS}, e^{-\phi}dV_{g_M\oplus g_{FS}}, m)$
satisfies the $AS_3$ condition.

\section{Appendix : Computing $F_3$}

In this Appendix, we prove (\ref{1.10})
by computing the term $F_3$ at a metric-measure structure $(g_\infty,\phi_\infty)$ with constant weighted scalar curvature.
First we will show that $F_1(v)=F_2(v)=0$.
To check this,
notice that $DF(w)[v]=DE(\Psi(w))[D\Psi(w)[v]]$.
Thus $DF(0)=0$ since $DE(1)=0$,
as $1$ is a critical point of the functional $E$ (by assumption, $(g_\infty,\phi_\infty)\in \mathcal{CWSC}$)
and $\Psi(0)=1$. Therefore, $F_1=0$.
Similarly,
$D^2F(w)[v,u]=D^2E(\Psi(w))[D\Psi(w)[u], D\Psi(w)[v]]
+\langle DE(\Psi(w)),D^2\Psi(w)[v,u] \rangle$.
When setting $w=0$, $\Psi(0)=1$, $D\Psi(0)=id$, and
\begin{equation*}
\begin{split}
D^2F(0)[v,u]&=
D^2E(1)[u,v]
+\langle DE(1),D^2\Psi(0)[v,u] \rangle\\
&=-\frac{8}{n+m-2}\langle \mathcal{L}_\infty u,v\rangle
+\langle DE(1),D^2\Psi(0)[v,u] \rangle.
\end{split}
\end{equation*}
As before, the second term vanishes. The first term
vanishes because $u$ is in the kernel of $\mathcal{L}_\infty$ by assumption. Therefore, $F_2=0$.

To compute $D^3F(0)$, we may in fact compute $D^3\tilde{F}(0)$,
where $\tilde{F}:\Lambda_0\to\mathbb{R}$ is defined by
$\tilde{F}(v)=E(1+v)$. We first compute $D^3F$:
\begin{equation}\label{4.2}
\begin{split}
D^3F(w)[v,u,z]&=D^3E(\Psi(w))[D\Psi(w)[v],D\Psi(w)[u],D\Psi(w)[z]]\\
&\quad+D^2E(\Psi(w))[D^2\Psi(w)[u,z],D\Psi(w)[v]]\\
&\quad+D^2E(\Psi(w))[D\Psi(w)[u],D^2\Psi(w)[v,z]]\\
&\quad+D^2E(\Psi(w))[D\Psi(w)[z],D^2\Psi(w)[v,u]]\\
&\quad+\langle DE(\Psi(w)), D^3\Psi(w)[v,u,z]\rangle.
\end{split}
\end{equation}
Setting $w=0$, and using similar considerations as before (in particular noting that
$D^2E(1)[\cdot]$ is self-adjoint), we obtain $D^3F(0)[v,u,z]=D^3E(1)[v,u,z]$.
Performing the same computation for $D^3\tilde{F}(0)$ yields the same result.
Next, we compute $D^3\tilde{F}(0)$.
In the rest of this section all integrals are taken with respect to $e^{-\phi_\infty}dV_{g_\infty}$. First we recall that
\begin{equation}\label{recall1}
  \frac{1}{2}DE(w)[v]=\int_M\left(a_{n,m}\Delta_{\phi_\infty}w+R_{\phi_\infty}^mw-r_\phi^mw^\frac{n+m+2}{n+m-2}\right)v
\end{equation}

where
\begin{equation*}
 a_{n,m}=-\frac{4(n+m-1)}{n+m-2},\quad  (g,\phi)=(w^\frac{4}{n+m-2}g_\infty,\phi_\infty-\frac{2m}{n+m-2}\ln w).
\end{equation*}
Because $r_\phi^m=E(w)$, eq \eqref{recall1} can be rewritten as
\begin{equation*}
  \frac{1}{2}DE(w)[v]=\int_M\left(a_{n,m}\Delta_{\phi_\infty}w+R_{\phi_\infty}^mw-E(w)w^\frac{n+m+2}{n+m-2}\right)v.
\end{equation*}
So the second differential of the functional $E$ can be computed as follows:

\begin{equation*}
  \begin{split}
    \frac{1}{2}D^2&E(u)[v,w]\\
    =&\frac{1}{2}\left.\frac{d}{dt}\right|_{t=0}DE(u+tw)[v]\\
    =&\left.\frac{d}{dt}\right|_{t=0}\int_M\left(a_{n,m}\Delta_{\phi_\infty}(u+tw)+R_{\phi_\infty}^m(u+tw)-E(u+tw)(u+tw)^\frac{n+m+2}{n+m-2}\right)v\\
    =&\int_M\left(a_{n,m}\Delta_{\phi_\infty}w+R_{\phi_\infty}^mw-DE(u)[w]u^\frac{n+m+2}{n+m-2}-\frac{n+m+2}{n+m-2}E(u)u^\frac{4}{n+m-2}w\right)v.
  \end{split}
\end{equation*}
Since $DE(1)=0$ and $E(1)=R_{\phi_\infty}^m$, we have
\begin{equation*}
  \frac{1}{2}D^2E(g_\infty,\phi_\infty)[v,w]=\int_M\left(-\frac{4(n+m-1)}{n+m-2}\Delta_{\phi_\infty}w-\frac{4}{n+m-2}w\right)v.
\end{equation*}
Similarly, the third differential of the functional $E$ can be computed as follows:
\begin{equation*}
  \begin{split}
    \frac{1}{2}D^3&E(u)[v,w,z]\\
    =&\frac{1}{2}\left.\frac{d}{dt}\right|_{t=0}D^2E(u+tz)[v,w]\\
    =&\left.\frac{d}{dt}\right|_{t=0}\int_M\left(a_{n,m}\Delta_{\phi_\infty}w+R_{\phi_\infty}^mw-DE(u+tz)[w](u+tz)^\frac{n+m+2}{n+m-2}\right.\\
      &\qquad\qquad\qquad \left.-\frac{n+m+2}{n+m-2}E(u+tz)(u+tz)^\frac{4}{n+m-2}w\right)v\\
    =&\int_M\left(-D^2E(u)[w,z]u^\frac{n+m+2}{n+m-2}-\frac{n+m+2}{n+m-2}DE(u)[w]u^\frac{4}{n+m-2}z\right.\\
    &\qquad \left.-\frac{n+m+2}{n+m-2}DE(u)[z]u^\frac{4}{n+m-2}w-\frac{4(n+m+2)}{(n+m-2)^2}E(u)u^{-\frac{n+m-6}{n+m-2}}zw\right)v.
  \end{split}
\end{equation*}
Therefore, $v,w,z\in\Lambda_0$, we have
\begin{equation*}
  D^3E(1)[v,w,z]=-\frac{8(n+m+2)}{(n+m-2)^2}R_{\phi_\infty}^m\int_Mvwz e^{-\phi_\infty}dV_{g_\infty}.
\end{equation*}

\section{Acknowledgement}

The first author was supported
by Basic Science Research Program through the National Research Foundation of Korea (NRF) funded by the Ministry of Education, Science and Technology (Grant No. 2020R1A6A1A03047877 and\\ 2019R1F1A1041021), and by Korea Institute for Advanced Study (KIAS) grant
funded by the Korea government (MSIP). The second author was supported by a KIAS Individual Grant (SP070701) via the Center for Mathematical Challenges at Korea Institute for Advanced Study.







\bibliographystyle{amsplain}

\end{document}